%% file: extended.tex
\documentclass[11pt]{article}

\usepackage[utf8]{inputenc}
\usepackage{amsmath}
\usepackage{amssymb}
\usepackage{amsthm}
\usepackage{graphicx}
\usepackage{float}
\usepackage{color}
\usepackage{bm}
\usepackage{mathtools}
\usepackage{geometry}
\usepackage{caption}
\usepackage{subcaption}
\usepackage{tikz}
\usepackage{blindtext}

\usepackage{hyperref}
\hypersetup{
	colorlinks=false,
	linkcolor=black,
	hyperindex=false,
	linktocpage=true,
}

\everymath{\displaystyle}
\newcommand{\pd}[2]{\displaystyle \frac{ \partial #1}{ \partial #2}}

\newcommand{\varpd}[2]{\displaystyle \frac{ \delta #1}{ \delta #2}}
\renewcommand{\div}{\mathrm{div}}

\DeclarePairedDelimiter\abs{\lvert}{\rvert}

\newcounter{num}
\newtheorem{remark}{Remark}[section]

\newcommand{\footremember}[2]{%
	\footnote{#2}
	\newcounter{#1}
	\setcounter{#1}{\value{footnote}}%
}
\newcommand{\footrecall}[1]{%
	\footnotemark[\value{#1}]%
}

\title{Extended Lagrangian approach for the numerical study of multidimensional dispersive waves: applications to the Serre-Green-Naghdi equations}
\author{%
	Sergey Tkachenko\footremember{imt}{Institut de Math\'{e}matiques de Toulouse, UMR CNRS 5219, INSA de Toulouse, 31077 Toulouse Cedex 4, France. Corresponding author, email: tkachenk@insa-toulouse.fr (sergeytkachenko012@gmail.com)}%
	\and Sergey Gavrilyuk\footremember{amu}{Aix-Marseille Universit\'e, CNRS UMR 7343 IUSTI, 5 rue Enrico Fermi, 13453 Marseille, France}%
	\and Jacques Massoni\footrecall{amu}%
}

\begin{document}

	\maketitle
	\input{sections/0_abstract}
	\section{Introduction}
		\input{sections/1_Introduction}
	
	\parindent=0cm
	\section{Dispersive models}\label{EXT_section_dispersive_models}
		\input{sections/2_Dispersive_models}

	\section{Extended Lagrangian formulation}\label{EXT_section_extended_lagrangian}
		\input{sections/3_Extended_Lagrangian_formulation}

	\section{Numerical resolution}\label{EXT_section_numerical_resolution}
		\input{sections/4_Numerical_resolution}
		
	\section{Numerical results}\label{EXT_section_results}
		\input{sections/5_Numerical_results}
	
	\section{Conclusion}
		\input{sections/6_Conclusion}
	
	\setcounter{section}{0}
	\renewcommand{\thesection}{\Alph{section}}
	\renewcommand{\theHsection}{2\Alph{section}}
	\input{sections/Appendix}

	\bibliographystyle{plain}
	\phantomsection
	\addcontentsline{toc}{chapter}{Bibliography}
	\bibliography{library}

\end{document}

%% file: sections/0_abstract.tex
\begin{abstract}
	
In this paper we study two multidimensional nonlinear dispersive systems: the Serre-Green-Naghdi (SGN) equations describing dispersive shallow water flows, and Iordanskii-Kogarko-Wijngaarden (IKW) equations describing fluids containing small compressible gas bubbles. These models are Euler-Lagrange equations for a given Lagrangian and share common mathematical structure, namely the dependence of the pressure on material derivatives of macroscopic variables. We develop a generic dispersive model such that SGN and IKW systems become its special cases if only one specifies the appropriate Lagrangian, and then use the extended Lagragian approach proposed in Favrie and Gavrilyuk (2017) to build its hyperbolic approximation. The new approximate model is unconditionally hyperbolic for both SGN and IKW cases, and accurately describes dispersive phenomena, which allows to impose discontinuous initial data and study dispersive shock waves. We consider the 2-D hyperbolic version of SGN system as an example for numerical simulations and apply a second order implicit-explicit scheme in order to numerically integrate the system. The obtained 1-D and 2-D results are in close agreement with available exact solutions and numerical tests.

\textbf{Keywords: } dispersive shallow water equations, bubbly fluids, Euler-Lagrange equations, hyperbolic conservation laws, multidimensional waves, implicit-explicit numerical methods

\end{abstract}

%% file: sections/1_Introduction.tex
A number of nonlinear dispersive systems possess a variational formulation, i.e. they are Euler-Lagrange equations for the Hamilton action which contains all physical information about the system. The dynamics of such a system depends only on the associated Lagrangian which is the difference between kinetic and potential energies. Being the most generic principle of mechanics, the Hamilton's principle states that the trajectory of a dissipationless system is a stationary point of the action functional (see \cite{Berdichevsky2009, Gavrilyuk2011}). In the dispersive case, the corresponding Lagrangian contains terms which depend not only on macroscopic variables of the system, but also on their spatial and temporal derivatives. A classical example is the Serre-Green-Naghdi equations which describe the propagation of long gravity waves on a surface of an inviscid irrotational incompressible fluid (Serre \cite{Serre1953}, Su and Gardner \cite{Su1969}, Green and Naghdi \cite{Green1974, Green1976}). The SGN system is a shallow water model of the second order of approximation with respect to the dimensionless small parameter $H_0 /L_0$, where $ L_0 $ and $ H_0 $ are characteristic horizontal and vertical lengths respectively. The associated average pressure depends not only on the water depth as in the classical Saint-Venant equations but also on the material derivatives of the water depth along the depth averaged velocity. These last terms take the acceleration of the free surface into account. Another model possessing the similar structure is the Iordanskii-Kogarko-Wijngaarden equations \cite{Iordanskii1960, Kogarko1961, Wijngaarden1968} describing fluids containing gas bubbles of small size. In this case, the pressure depends not only on the gas density but also on the material derivatives of the density up to second order (this dynamic equation is called Rayleigh-Lamb equation). Analogous mathematical models also appear in the description of  shock wave propagation  in ductile porous metals where the micro-inertia effects related with the rapid variation of the porosity become important (see \cite{Molinari2017}).

Since these models share a common mathematical structure and describe qualitatively similar phenomena, our numerical simulations will concern the SGN equations as the most studied example of such systems. Le M\'{e}tayer et al. \cite{LeMetayer2010} developed a hybrid finite volume/finite difference scheme based on the conservative formulation of the SGN model. The hyperbolic part of the model is treated by a Godunov type method, and the dispersive part is treated by a finite difference scheme. At each time step, after the resolution of the hyperbolic part, an elliptic operator is inverted in the whole numerical domain. This strategy was also used by Bonneton et al. \cite{Bonneton2011a} for the case of varying topography. They introduced two high-order numerical approaches to the numerical resolution of the SGN model. The first one is a high order hybrid finite volume/finite difference method based on the splitting scheme mentioned above, and second one is based on the quasi-conservative form of SGN equations. Both of them include a special way to handle the wave breaking: at some point, when the wave slope becomes critical, the wave is ``ready to break'', and the numerical model switches from the SGN equations to first-order hyperbolic shallow water equations which admit shock waves. The switch is performed locally in space and time when the energy dissipation is high. In another article, Bonneton et al. \cite{Bonneton2011} proposed a similar high-order approach based on the finite volume/finite difference splitting and also introduced a formulation with an improved dispersion relation. Li et al. \cite{Li2014} developed a similar hybrid approach which includes the inversion of a global operator but instead of finite volumes/finite differences they used the continuous Galerkin/finite element hybrid scheme. Chazel et al. \cite{Chazel2011} introduced a three-parameter model which tends to original SGN equations in the long-wave limit and which dispersive relation is very close to the one of the Euler equations for the full water wave problem, with a proper choice of parameters. The numerical method is also based on the high-order hybrid approach mentioned above. There are also 2-D extensions of those numerical approaches, namely the 2-D approach of Lannes and Marche \cite{Lannes2015}. Although the mentioned approaches accurately describe the dispersive phenomena, the elliptic parts of the schemes need to be treated globally, which affects the numerical performance of the methods.

The first attempt to replace the fully nonlinear 1-D second-order models of the shallow-water theory by hyperbolic approximations was made by Liapidevskii and Gavrilova \cite{Liapidevskii2008}. They proposed a conditionally hyperbolic model using a relaxation technique. Favrie and Gavrilyuk \cite{Favrie2017} developed a new approach for the ``hyperbolization'' based on the variational structure. It consists in modifing the original Lagrangian (``master'' Lagrangian) by introducing a one-parameter family of new extended Lagrangians. This Lagrangian contains new ``penalized'' macroscopic variables : these new variables tend to old variables in some limit. Here the variational formulation becomes extremely useful, since one only needs to modify the Lagrangian and then apply the Hamilton's principle to the corresponding action. The new governing equations are unconditionally hyperbolic. The mathematical justification of the ``penalisation'' technique was given by  Duch\^{e}ne \cite{Duchene2019}. This approach was as well successfully applied to the 1-D nonlinear Schr\"{o}dinger equation by Dhaouadi et al. \cite{Dhaouadi2018}. Let us specifically mention the approximate hyperbolic systems including the varying topography, namely the scale invariant relaxation model by Guermond et al. \cite{Guermond2019} and a quasi-incompressible model by Richard \cite{Richard2021}. Recently Dumbser et al. \cite{Gavrilyuk2021} introduced the implementation of high order ADER discontinuous Galerkin schemes for the 2-D extension of the model \cite{Favrie2017} with varying topography. 

In this paper we develop a multidimensional dispersive model which unifies both SGN and IKW bubbly fluids systems under one single formulation. Then, we apply the the extended Lagrangian method \cite{Favrie2017} in order to build its extended hyperbolic extension. For numerical simulations, we consider a particular case of the new extended system corresponding to the 2-D flat bottom case of the extended SGN system \cite{Favrie2017}. In order to integrate the system numerically, we make use of the second-order implicit-explicit approach (IMEX). Originally developed in \cite{Ascher1997, Pareschi2001, Pareschi2005}, it demonstrated robustness and high precision for hyperbolized one-dimensional dispersive systems with stiff source terms \cite{Dhaouadi2020, Tkachenko2020a, Richard2021}. In this article, we extend this approach for 2-D simulations.

The article is organized as follows. In section \ref{EXT_section_dispersive_models} we introduce the generic dispersive model, unifying the SGN and IKW systems; its hyperbolic version is presented in section \ref{EXT_section_extended_lagrangian}. The numerical methods are described in section \ref{EXT_section_numerical_resolution}, followed by section \ref{EXT_section_results}, describing the numerical results.

%% file: sections/2_Dispersive_models.tex
\subsection{Serre-Green-Naghdi model}

	The SGN model \cite{Green1974, Green1976, Serre1953, Su1969} is given by the following equations:
	\begin{equation}\label{1-disp_sgn-equations}
		\begin{aligned}
		&\pd{h}{t} + \div(h \textbf{u}) = 0,\\
		&\pd{h \textbf{u}}{t} + \div(h \textbf{u} \otimes \textbf{u} + p I) = 0, \qquad p = \frac{gh^2}{2} + \frac{1}{3}h^2\ddot{h}.
		\end{aligned}
	\end{equation}
		Here $ h $ is the water depth, $ \textbf{u} = (u_1, u_2)^T $ is the horizontal velocity averaged over the water depth, $ g $ is the gravity acceleration and dots denote material derivatives:
	\begin{equation*}
		\dot{h} = \pd{h}{t} + \textbf{u} \cdot \nabla h, \qquad \ddot{h} = \pd{\dot{h}}{t} + \textbf{u} \cdot \nabla \dot{h}.
	\end{equation*}
	The pressure is non-hydrostatic and depends not only on macroscopic variable $ h $ but also on its second-order material derivative. The momentum equation in \eqref{1-disp_sgn-equations} is Euler-Lagrange equation for the Lagrangian (see \cite{Gavrilyuk2011}:
	\begin{equation}\label{1-disp_sgn-lagrangian}
		\mathcal{L} = \int\limits_{D(t)}\left(\frac{h\!\left|\textbf{u}\right|^2}{2} - W(h,\dot{h})\right)dD,
	\end{equation}
	where the potential $ W $ is given by:
	\begin{equation}\label{1-disp_sgn-potential}
		W(h,\dot{h}) = \frac{gh^2}{2} - \frac{h\dot{h}^2}{6}.
	\end{equation}
	The SGN system \eqref{1-disp_sgn-equations} admits the energy conservation law:
	\begin{equation}\label{1-disp_sgn-energy-conservation}
		\pd{ E}{t} + \div \left( E\textbf{u} + p \textbf{u} \right) = 0,
	\end{equation}
	where the total energy of the system is given by:
	\begin{equation}\label{1-disp_sgn-energy}
	E = \frac{h\!\left|\textbf{u}\right|^2}{2} + \frac{h\dot{h}^2}{6} + \frac{gh^2}{2}.
	\end{equation}

\subsection{Iordanskii-Kogarko-Wijngaarden model}

	Consider an incompressible fluid of density $ \rho_{10} = const $ containing bubbles of compressible gas. We will consider the IKW system \cite{Iordanskii1960, Kogarko1961, Wijngaarden1968} respecting the following assumptions. First, all bubbles have the same radius $ R $ at a given point of space. Second, surface tension, viscosity and heat conduction are neglected. Third, the bubble radius is significantly smaller than the inter-bubble distance $ d $ which, in turn, is significantly smaller than the scale of motion $ l $:
	\begin{equation*}
		R \ll d \ll l.
	\end{equation*}
	The governing equations are as follows:
	\begin{equation}\label{1-disp_ikw-equations}
		\begin{aligned}
			&\pd{\rho}{t} + \div(\rho \textbf{u}) = 0,\\
			&\pd{\rho \textbf{u}}{t} + \div(\rho \textbf{u} \otimes \textbf{u} + p I) = 0,\\
			& p = p_2 + \rho_{10}\left( \frac{3}{2} \dot{R}^2 + R \ddot{R} \right),\\
			& \pd{N}{t} + \div (N\textbf{u}) = 0.
		\end{aligned}
	\end{equation}
	Here $ \textbf{u} = (u_1,u_2,u_3)^T $ is the mean velocity of the mixture motion, $ \rho = \alpha_2 \rho_2 + \alpha_1 \rho_{10} $ is the mixture density, where $ \rho_2 $ is the gas density, and $ \alpha_1$ and $\alpha_2 $ are the volume fractions of liquid and gas correspondingly, such that $ \alpha_1 + \alpha_2 = 1 $. The gas pressure inside a bubble is denoted by $ p_2 $, $ R $ is the bubble radius, and $ N $ is the number of bubbles per unit volume. The volume fraction and density of gas are expressed as follows:
	\begin{equation*}
		\alpha_2 = \frac{4 \pi R^3}{3} N, \qquad \rho_2 = \frac{Y_2\rho}{\alpha_2} = \frac{Y_2}{\frac{\alpha_2}{\rho}} = \frac{Y_2}{\frac{4}{3}\pi R^3 n}.
	\end{equation*}
	Then, the system \eqref{1-disp_ikw-equations} is closed if $ n = const $ (the bubbles neither disappear nor appear), $ Y_2 = const $ (there is no mass exchange between phases). Here we introduced the mass fractions
	\begin{equation*}
		Y_1 = \frac{\alpha_1 \rho_{10}}{\rho}, \qquad Y_2 = \frac{\alpha_2 \rho_2}{\rho}, \qquad Y_1 + Y_2 = 1,
	\end{equation*}
	and	the number $ n $ of bubbles per unit mass:
	\begin{equation*}
		n = \frac{N}{\rho}.
	\end{equation*}
	Then, using the identity $ \alpha_1 + \alpha_2 = 1 $, we link the mixture density $ \rho $ to $ R $ as follows:
	\begin{equation}\label{1-disp_radius-density}
		\frac{4}{3} \pi R^3 = \frac{1}{n} \left( \frac{1}{\rho} - \frac{Y_1}{\rho_{10}} \right),
	\end{equation}
	thus the bubble radius can be expressed as a function of density. We introduce the micro-inertial kinetic energy of the fluid, appearing due to oscillations of $ N $ bubbles \cite{Iordanskii1960, Kogarko1961, Wijngaarden1968}:
	\begin{equation}\label{1-disp_bubbles-micro-energy}
		2 \pi R^3 N \rho_{10}\dot{R}^2.
	\end{equation}
	The momentum equation and pressure equation in  \eqref{1-disp_ikw-equations} are Euler-Lagrange equations for the Lagrangian
	\begin{equation}\label{1-disp_ikw-lagrangian}
		\mathcal{L} = \int\limits_{D(t)} \left(\frac{ \rho\left| \textbf{u} \right|^2 }{ 2 } + 2 \pi R^3 N \rho_{10}\dot{R}^2 - \rho Y_2 \varepsilon_2 (\rho_2) \right)d D.
	\end{equation}
	Here $ \varepsilon_2 $ is the specific energy of the gas. Now we want to express the Lagrangian as a function of $ \textbf{u} $, $ \rho $ and $ \dot{\rho} $. First of all, we will express $ \rho Y_2 \varepsilon_2 (\rho_2) $ it in terms of $ \rho $. We will suppose that the gas is polytropic:
	\begin{equation*}
		p_2 = p_0 \left( \frac{ V_0 }{ V_2 } \right)^{\gamma}.
	\end{equation*}
	Here $ V_2 = \frac{ 4 }{ 3 } \pi R^3 $ is a single bubble volume, $ V_0 = \frac{ 4 }{ 3 } \pi R_0^3 $ is the initial bubble volume and $ \gamma > 1 $. The isentropic Gibbs identity written in volume units reduces to:
	\begin{equation*}
		d E_2 + p_2 d V_2 = 0,
	\end{equation*}
	where $ E_2 $ is the gas volume energy. Integration of the Gibbs identity over the volume occupied by gas gives:
	\begin{equation*}
	E_2 = p_0 \int \left(\frac{V_0}{V_2}\right)^{-\gamma} d V_2 = \frac{ p_0 V_0^{\gamma} }{ \gamma - 1 } V_2^{ - \gamma + 1 } = \frac{ p_2 V_2 }{ \gamma - 1 }.
	\end{equation*}
	Or, written in specific quantities with $ \tau_2 = 1/\rho_2 $, the gas specific energy $ \varepsilon $ reads:
	\begin{equation*}
	\varepsilon_2 = \frac{ p_2 \tau_2 }{ \gamma - 1 }.
	\end{equation*}
	Finally, we express the potential energy of the gas fraction as follows:
	\begin{equation*}
		\rho Y_2 \varepsilon_2 = \alpha_2 \rho_2 \varepsilon_2 = \frac{ \alpha_2 p_2 }{ \gamma - 1 } = \frac{ \frac{4}{3} \pi R^3 N }{ \gamma - 1 } p_0 \left( \frac{ V_0 }{ V_2 } \right)^{3 \gamma} = \frac{ \frac{4}{3} \pi R^3 \rho n }{ \gamma - 1 } p_0 \left( \frac{ R_0 }{ R } \right)^{3 \gamma}.
	\end{equation*}

	Since the bubble radius $ R $ depends on $ \rho $ via \eqref{1-disp_radius-density}, we can rewrite the Lagrangian \eqref{1-disp_ikw-lagrangian} with the potential $ W $ taken as a function of $ \rho $ and $ \dot{\rho} $:
	\begin{equation*}
		\mathcal{L} = \int\limits_{D(t)} \Bigg(  \frac{\rho\!\left|\textbf{u}\right|^2}{2} - W(\rho,\dot{\rho}) \Bigg)~d D,
	\end{equation*}
	where $ W(\rho, \dot{\rho}) $ reads:
	\begin{equation}\label{1-disp_ikw-potential}
		W(\rho,\dot{\rho}) = \rho \left(\frac{ \frac{4}{3} \pi R^3 n }{ \gamma - 1 } p_0 \left( \frac{ R_0 }{ R } \right)^{3 \gamma} - 2 \pi n \rho_{10} R^3\dot{R}^2\right).
	\end{equation}
	Let us notice that such formulation of the Lagrangian is completely analogous to \eqref{1-disp_sgn-lagrangian}. IKW model also admits the energy conservation law in the same form as \eqref{1-disp_sgn-energy-conservation}, with total energy $ E $ given by:
	\begin{equation}\label{1-disp_ikw-total-energy}
		E = \rho \left(\frac{\left|\textbf{u}\right|^2}{2} + 2 \pi n \rho_{10} R^3\dot{R}^2 + \frac{ \frac{4}{3} \pi R^3 n }{ \gamma - 1 } p_0 \left( \frac{ R_0 }{ R } \right)^{3 \gamma}\right).
	\end{equation}
	Summarizing this and the previous sections, we can say that both SGN and IKW models are Euler-Lagrange equations with Lagrangians which possess the same structure, namely the dependency on the potential on macroscopic variables and its material derivatives of the first order. This allows us to write both models \eqref{1-disp_sgn-equations} and \eqref{1-disp_ikw-equations} using one generic formulation:
	\begin{equation}\label{1-disp_generic-formulation}
		\begin{aligned}
			&\pd{\rho}{t} + \div(\rho \textbf{u}) = 0,\\
			&\pd{\rho \textbf{u}}{t} + \div(\rho \textbf{u} \otimes \textbf{u} + p I) = 0, \qquad p = \rho\varpd{W}{\rho} - W.
		\end{aligned}		
	\end{equation}
	Where $ \varpd{W}{\rho} $ is the variational derivative of $ W $ \cite{Gavrilyuk2004, Gavrilyuk2011}:
		\begin{equation*}
			\varpd{W}{\rho} = \pd{W}{\rho} - \pd{}{t}\left(\pd{W}{\dot{\rho}}\right) - \div\left(\pd{W}{\dot{\rho}} \textbf{u}\right).
		\end{equation*}
	If $ \textbf{u} = (u_1,u_2)^T $, $ \rho = h $ and the potential $ W $ is defined by \eqref{1-disp_sgn-potential}, then \eqref{1-disp_generic-formulation} is the Serre-Green-Naghdi model \eqref{1-disp_sgn-equations}. If $ \textbf{u} = (u_1,u_2,u_3)^T $ and $ W $ is given by \eqref{1-disp_ikw-potential}, then the equations become the Iordanskii-Kogarko-Wijngaarden model \eqref{1-disp_ikw-equations}.

%% file: sections/3_Extended_Lagrangian_formulation.tex
The original idea of authors \cite{Favrie2017} was to introduce a new non-equilibrium variable $ \eta $ for the SGN model which tends to the fluid depth $ h $ in some limit and replaces it in the micro-inertial kinetic energy term $ \frac{h\dot{h}^2}{6} $. Then, a penalization term with a large parameter $ \lambda \gg 1 $ is added to the new extended Lagrangian to assure this convergence:

\begin{equation*}
\tilde{L}(\textbf{u}, h,\eta,\dot{\eta}) = \int\limits_{D(t)} \Bigg( \frac{h |\textbf{u}|^2}{2} + \frac{h\dot{\eta}^2}{6} - \frac{gh^2}{2} - \frac{\lambda h}{6}\left( \frac{\eta}{h} -1 \right)^2 \Bigg) d D.
\end{equation*}
This chapter is aimed to propose a hyperbolic extension for the generic system \eqref{1-disp_generic-formulation} using the extended Lagrangian approach. We first put the bubbly fluids micro-inertial energy term \eqref{1-disp_bubbles-micro-energy} in the same quadratic form as it is done for the SGN model:

\begin{equation}\label{2-extl_Q-equation}
	2 \pi n \rho_{10} \rho R^3\dot{R}^2 = \frac{1}{2}\beta\rho \,\dot{Q}^2,
\end{equation}
Here we introduce a new function $ Q $ which depends explicitly on $ \rho $ via \eqref{1-disp_radius-density}:
\begin{equation}\label{2-extl_Q-definition}
	Q(\rho) = \frac{2}{5} R^{\frac{5}{2}}(\rho), \qquad \dot{Q} = \frac{d Q}{d \rho} \dot{\rho}.
\end{equation}
The definition for $ \beta $ not depending on $ \rho $ comes automatically:
\begin{equation}\label{2-extl_ikw-beta-definition}
	\beta = 4 \pi n \rho_{10}.
\end{equation}
Thus, we can now rewrite the generic potential $ W $ from \eqref{1-disp_generic-formulation} into the following form:

\begin{equation}\label{2-extl_generic-potential}
	W(\rho,\dot{\rho}) = \rho \left(\varepsilon(\rho) - \frac{1}{2}\beta\dot{Q}^2\right), \qquad Q = Q(\rho).
\end{equation}
We will exploit the notation of $ \rho $ for both Serre-Green-Naghdi and Iordanskii-Kogarko-Wijngaarden models, supposing that when we speak about the former, $ \rho $ stands for the water depth $ h $. If we take the specific energy $ \varepsilon(\rho) $, the function $ Q(\rho) $ and the constant $ \beta $ as follows:
\begin{equation}\label{2-extl_SGN-functions-definitions}
	\varepsilon(\rho) = \frac{g\rho}{2}, \qquad Q(\rho) = \rho, \qquad \beta = \frac{1}{3},
\end{equation}
then \eqref{2-extl_generic-potential} defines the potential of the Serre-Green-Naghdi model \eqref{1-disp_sgn-potential}. In the same way, expressions \eqref{2-extl_Q-definition}, \eqref{2-extl_ikw-beta-definition} define the bubbly fluids model potential \eqref{1-disp_ikw-potential} with the internal energy coming without any transformations:
\begin{equation}\label{2-extl_ikw-specific-energy}
	\varepsilon(\rho) = \frac{ \frac{4}{3} \pi R^3(\rho) n }{ \gamma - 1 } p_0 \left( \frac{ R_0 }{ R(\rho) } \right)^{3 \gamma}.
\end{equation}
For the sake of simplicity we suppose that $ R_0 $, $ p_0 $, $ n $ and $ Y_l $ are identically constant. Now that the potentials of both SGN and IKW models are written under one generic form, we will employ the extended Lagrangian approach \cite{Favrie2017} to construct a new model, approximating the generic formulation \eqref{1-disp_generic-formulation}. We replace $ Q(\rho) $ in \eqref{2-extl_generic-potential} by the new variable $ \eta $:
\begin{equation}\label{2-extl_eta-tends-to-Q}
	\eta \rightarrow Q(\rho),
\end{equation}
or, equally:
\begin{equation*}
	f(\eta) \rightarrow \rho,
\end{equation*}
where $ f = Q^{-1} $. In order to provide this convergence we add a generic penalty term
\begin{equation}\label{2-extl_penalty-term-definition}
 \frac{a\lambda \rho }{2} \left( \frac{ f(\eta) }{ \rho } - 1 \right)^2, \qquad a = const,
\end{equation}
to the potential, also replacing $ \dot{Q}(\rho) $ by $ \dot{\eta} $, so that the new extended potential $ \tilde{W}(\rho,\eta,\dot{\eta}) $ becomes:
\begin{equation*}
	\tilde{W}(\rho,\eta,\dot{\eta}) =  \rho\varepsilon(\rho) - \frac{\beta\rho\dot{\eta}^2}{2} + \frac{a\lambda \rho }{2} \left( \frac{ f(\eta) }{ \rho } - 1 \right)^2.
\end{equation*}
If $ \lambda $ goes to infinity, then $ \left( \frac{ f(\eta) }{ \rho } - 1 \right)^2 $ tends to zero, which automatically provides the convergence \eqref{2-extl_eta-tends-to-Q}. We define the generic extended Lagrangian $ \tilde{L}(\textbf{u},\rho,\eta,\dot{\eta}) $ as follows:
\begin{equation*}
 \tilde{L}(\textbf{u},\rho,\eta,\dot{\eta}) = \int\limits_{D(t)} \Bigg( \frac{\rho |\textbf{u}|^2}{2} - \tilde{W}(\rho,\eta,\dot{\eta}) \Bigg) d D.
\end{equation*}
The corresponding Euler-Lagrange equations for the extended Lagrangian are:
\begin{equation}\label{2-extl_extended-system-eta-ddot}
	\begin{aligned}
		&\pd{\rho}{t} + \div(\rho \textbf{u}) = 0,\\
		& \pd{\textbf{u}}{t} + \pd{\textbf{u}}{\textbf{x}}\textbf{u} + \frac{1}{\rho} \nabla p = 0,\\
		& \ddot{\eta} = - \frac{a\lambda f'(\eta)}{\beta\rho} \left(\frac{f(\eta)}{\rho} - 1 \right),\\
	\end{aligned}
\end{equation}
where $ f'(\eta) = \frac{df(\eta)}{d\eta} $. The new ``pressure'' $ p $ depends only on $ \rho $ and $ \eta $:
\begin{equation*}
	p(\rho,\eta) = \rho \pd{\tilde{W}}{\rho}(\rho,\eta,\dot{\eta}) - \tilde{W}(\rho,\eta,\dot{\eta}) = \rho^2 \varepsilon'(\rho) - a\lambda f(\eta)\left(\frac{f(\eta)}{\rho} - 1\right),
\end{equation*}
which means that the fluxes do not depend on derivatives anymore. The third equation in \eqref{2-extl_extended-system-eta-ddot} is of the second order, so we introduce a new variable $ w = \dot{\eta} = \pd{\eta}{t} + (\textbf{u}\cdot\nabla)\eta $ to rewrite \eqref{2-extl_extended-system-eta-ddot} as a first order system:

\begin{equation}\label{2-extl_extended-system}
\begin{aligned}
&\pd{\rho}{t} + \div(\rho \textbf{u}) = 0,\\
& \pd{\textbf{u}}{t} +\pd{\textbf{u}}{\textbf{x}}\textbf{u} + \frac{1}{\rho} \nabla p = 0, \qquad p = \rho^2\varepsilon'(\rho) - a\lambda f(\eta)\left(\frac{f(\eta)}{\rho} - 1\right),\\
& \pd{\eta}{t} + (\textbf{u}\cdot\nabla)\eta = w,\\
& \pd{w}{t} + (\textbf{u}\cdot\nabla)w =  - \frac{a\lambda f'(\eta)}{\beta\rho} \left(\frac{f(\eta)}{\rho} - 1 \right).
\end{aligned}
\end{equation}
The eigenvalues of a 1-D system are given by:
\begin{equation*}
\mu_{1,2,3,4} = u_1, \qquad \mu_{5, 6} = u_1 \pm \sqrt{p_\rho}.
\end{equation*}
The full 3-D system is hyperbolic if the local ``sound speed'' is positive:
\begin{equation*}
\pd{p}{\rho} = \rho\big(\rho\varepsilon(\rho)\big)'' + a\lambda \frac{f^2(\eta)}{\rho^2} > 0.
\end{equation*}
Contact characteristics $ \mu_{1, 2, 3, 4} = u_1 $ are \textit{linearly degenerate}:
\begin{equation*}
	\mathbf{r}_k\cdot \nabla_{\mathbf{U}} (\mu_k)  \equiv 0, \qquad k = 1,2,3,4, \qquad \nabla_{\textbf{U}} = \left(\partial_\rho, \partial_{u_1}, \partial_{u_2}, \partial_{u_3}, \partial_\eta, \partial_w\right)^T.
\end{equation*}
``Sound'' characteristics $ \mu_{5, 6} = u_1 \pm \sqrt{p_\rho} $ are \textit{genuinely non-linear in the sense of Lax}\cite{Lax2006}:
\begin{equation*}
	\mathbf{r}_{5, 6}\cdot \nabla_{\mathbf{U}} \mu_{5, 6} = \pm \frac{ 1 }{ 2 } \frac{ \sqrt{p_{\rho}} }{ \rho } \left( 2 + \frac{ \rho p_{\rho\rho} }{ p_{\rho} } \right) \ne 0.
\end{equation*}
Here $ \textbf{r}_k $ are right eigenvectors of the 1-D system. For full hyperbolicity study see Appendix \ref{appendix:extended-model_hyperbolicity}. The system also admits a general energy conservation law (see Appendix \ref{appendix:extended-model_energy-conservation}):
\begin{equation*}
	\pd{ E}{t} + \div \left( E\textbf{u} + p \textbf{u} \right) = 0,
\end{equation*}
where
\begin{equation*}
	E =  \frac{\rho\!\left|\textbf{u}\right|^2}{2} +  \frac{\beta\rho\dot{\eta}^2}{2} + \rho\varepsilon(\rho) + \frac{ a\lambda \rho }{ 2 } \left( \frac{ f(\eta) }{ \rho } - 1 \right)^2.
\end{equation*}
Throughout this work will call the system \eqref{2-extl_extended-system} with $ \textbf{u} = (u_1,u_2)^T $, $ \rho = h $, $ f(\eta) = \eta $, $ a = \frac{1}{3} $ and definitions from \eqref{2-extl_SGN-functions-definitions} the \textit{extended Serre-Green-Naghdi system} \cite{Favrie2017}:
\begin{equation}\label{2-extl_favrie-gavrilyuk-system}
	\begin{aligned}
		& \pd{h}{t} + \div(h\textbf{u})_x = 0,\\
		& \pd{\textbf{u}}{t} + \pd{\textbf{u}}{\textbf{x}}\textbf{u} + \frac{1}{h}\nabla\left(\frac{gh^2}{2} - \frac{\lambda\eta}{3}\left(\frac{\eta}{h}-1\right)\right) = 0,\\
		& \pd{\eta}{t} + (\textbf{u}\cdot\nabla)\eta = w,\\
		& \pd{w}{t} + (\textbf{u}\cdot\nabla)w = - \frac{\lambda}{h}\left(\frac{\eta}{h} - 1\right).
	\end{aligned}
\end{equation}
The system with $ \textbf{u} = (u_1,u_2,u_3) $, $ a = 1 $, $ \beta $ from \eqref{2-extl_ikw-beta-definition}, and $ \varepsilon(\rho) $ from \eqref{2-extl_ikw-specific-energy} will be called the \textit{extended Iordanskii-Kogarko-Wijngaarden system}:
\begin{equation}\label{2-extl_extended-ikw-system}
	\begin{aligned}
		&\pd{\rho}{t} + \div(\rho \textbf{u}) = 0,\\
		& \pd{\textbf{u}}{t} +\pd{\textbf{u}}{\textbf{x}}\textbf{u} + \frac{1}{\rho} \nabla \left(p_0\left(\frac{R_0}{R(\rho)}\right)^{3\gamma} - \lambda f(\eta)\left(\frac{f(\eta)}{\rho}-1\right)\right) = 0,\\
		& \pd{\eta}{t} + (\textbf{u}\cdot\nabla)\eta = w,\\
		& \pd{w}{t} + (\textbf{u}\cdot\nabla)w =  - \frac{\lambda f'(\eta)}{4 \pi n \rho_{10}\rho} \left(\frac{f(\eta)}{\rho} - 1 \right).
	\end{aligned}
\end{equation}
where $ f(\eta) = Q^{-1}(\eta) $ from \eqref{2-extl_Q-definition}:
\begin{equation*}
	f(\eta) = Q^{-1}(\eta) =\frac{1}{\frac{ Y_1 }{ \rho_{10} } + \frac{ 4 \pi n }{ 3 } \left( \frac{ 5 }{ 2 } \eta \right)^{\frac{ 6 }{ 5 }}}.
\end{equation*}

Both extended SGN \eqref{2-extl_favrie-gavrilyuk-system} and extended IKW \eqref{2-extl_extended-ikw-system} systems are unconditionally hyperbolic, since the following conditions are satisfied (see Appendix \ref{appendix:extended-model_hyperbolicity} for details):
\begin{equation*}
\begin{aligned}
& \pd{p_{ikw}}{\rho} = \frac{3 \gamma p_0}{4 \pi n \rho^2 R^3}\left(\frac{R_0}{R}\right)^{3\gamma} + \lambda \frac{f^2(\eta)}{\rho^2} > 0,\\
& \pd{p_{sgn}}{\rho} = g\rho + \frac{\lambda}{3}\frac{\eta^2}{\rho^2} > 0.
\end{aligned}
\end{equation*}

\begin{remark}
	The potential $ W(\rho,\dot{\rho}) $ can also be extended to $ \tilde{W} = \tilde{W}(\rho,c,\eta,\dot{\eta}) $, for some Lagrangian variable $ c(t,\textbf{x}) $ conserved along the trajectories: i.e. such that $ \dot{c} = 0 $. It would permit us to consider non-homogeneous media. For instance, it could be the initial space-dependent bubble radius $ R_0 $, number of bubbles per unit mass $ n $, initial pressure $ p_0 $, etc. Although they might not be identically constant, the Euler-Lagrange equations will stay the same, and one will only need to add new transport equations for these variables to obtain the full system.
\end{remark}

\begin{remark}
	The idea of the penalization technique is intuitively quite understandable, but its mathematical justification is not at all obvious. Such a justification for a hyperbolized SGN system was recently done by V. Duch\^{e}ne \cite{Duchene2019}.
\end{remark}

%% file: sections/4_Numerical_resolution.tex
Both models described above possess the same mathematical structure. However, in terms of visualization of physical processes, surface waves motion appears to be more intuitive since the evolution of fluid parameters is clearly observed, while those of bubbly fluids can only be measured. Hence, we will focus on numerical resolution of the extended SGN model. We rewrite \eqref{2-extl_favrie-gavrilyuk-system} in a conservative form, using the notation $ u_1 = u $, $ u_2 = v $:
\begin{equation}\label{3-num_ext-sgn-2d-conservative-vector}
\pd{\textbf{U}}{t} + \pd{\textbf{F}(\textbf{U})}{x} + \pd{\textbf{G}(\textbf{U})}{y} = \textbf{S}(\textbf{U}).
\end{equation}
where $ \textbf{U} $ is the vector of conservative variables, $ \textbf{S}(\textbf{U}) $ is the vector of source terms
\begin{equation}\label{3-num_ext-sgn-2d-conservative-vector_sources}
\textbf{U} = (h, hu, hv, h\eta, hw)^T, \qquad \textbf{S}(\textbf{U}) = \left(0, 0, 0, hw, - \lambda \left( \frac{\eta}{h} - 1\right) \right)^T,
\end{equation}
and $ \textbf{F}(\textbf{U}) $, $ \textbf{G}(\textbf{U}) $ are the flux vectors independently separated in $ x $ and $ y $ directions:
\begin{equation}\label{3-num_ext-sgn-2d-conservative-vector_fluxes-pressure}
\textbf{F} = (hu, hu^2 + p, huv, hu\eta, huw)^T, \qquad \textbf{G} = (hv, hvu, hv^2 + p, hv\eta, hvw)^T, \qquad p = \frac{gh^2}{2}-\frac{\lambda \eta}{3}\left(\frac{\eta}{h}-1\right).
\end{equation}

Consider a rectangular domain $ [x_l, x_r] \times [y_l, u_r]$. We divide the it into $ N_x \times N_y $ equal rectangular cells $ \left(C_{ij}\right)_{1<i<N_x, 1<j<N_y} $ with sides $ \Delta x = \frac{x_r - x_l}{N_x} $, $ \Delta y = \frac{y_r - y_l}{N_y} $, such that $ C_{i,j} =\left [x_{i-\frac{1}{2}},x_{i+\frac{1}{2}}\right]\times\left[y_{j-\frac{1}{2}},y_{j+\frac{1}{2}}\right] $, where $ x_{i\pm\frac{1}{2}} $ and $ y_{j\pm \frac{1}{2}} $ are the cell boundaries of the nodes of a regular Cartesian mesh:
\begin{equation*}
\left\{\begin{array}{ll}
x_i = (i/2 - 1)\Delta x, \qquad & 1<i<N_x,\\
y_j = (j/2 - 1)\Delta y, \qquad & 1<j<N_y.
\end{array}\right.
\end{equation*}
A single time step will be denoted $ \Delta t $, and the nodal value of any function $ u $ in the cell $ C_{i,j} $ at the moment $ t_n $ is denoted $ u_{ij}^n $. In the present article will compare two finite volume methods of first and second order correspondingly.
	
\subsection{First-order splitting}
	The first-order splitting method for the 1-D extended SGN system is introduced in the original work \cite{Favrie2017}, and we present it's straightforward 2-D extension:
	\begin{equation}\label{3-num_first-order-scheme}
		\begin{aligned}
			& \textbf{U}_{ij}^{(1)} = \textbf{U}_{ij}^n  - \Delta t\, \mathcal{G}\left( \textbf{U}_{ij}^n \right),\\
			& \textbf{U}_{ij}^{(2)} = ODE_{exact} \left( \textbf{U}_{ij}^n \right)\\
			& \textbf{U}_{ij}^{n+1} = \textbf{U}_{ij}^{(1)}  + \Delta t \, \textbf{S} \left(\textbf{U}_{ij}^{(2)} \right)
		\end{aligned}
	\end{equation}
	The first step resolves the homogeneous part of the system. Here $ \mathcal{G} $ is a 2-D operator taken as described in \cite{Toro2009}, i.e. the numerical solution is updated in $ x $ and $ y $ directions simultaneously in a single time step:
	\begin{equation}\label{3-num_hyperbolic operator}
		\mathcal{G}(\textbf{U}_{ij}^n) = \frac{1}{\Delta x} \left( \textbf{F}_{i+\frac{1}{2},j}^n - \textbf{F}_{i-\frac{1}{2},j}^n \right) + \frac{1}{\Delta y} \left( \textbf{G}_{i, j+\frac{1}{2}}^n - \textbf{G}_{i, j-\frac{1}{2}}^n \right).
	\end{equation}
	The intercell numerical fluxes $ \textbf{F}_{i+\frac{1}{2}, j}^n $ and $ \textbf{G}_{i, j+\frac{1}{2}}^n $ are obtained via resolution of the Riemann problem on the cell boundaries $ i \pm 1/2 $ and $ j \pm 1/2 $:
	\begin{equation*}
		\begin{aligned}
			& \textbf{F}_{i+\frac{1}{2}, j}^n = RP(\textbf{U}_{ij}^n,\textbf{U}_{i+1,j}^n),\\
			& \textbf{G}_{i, j+\frac{1}{2}}^n = RP(\textbf{U}_{ij}^n,\textbf{U}_{i,j+1}^n).
		\end{aligned}
	\end{equation*}
	We will utilize two Riemann solvers in this article, the choice depends on the problem to consider. Rusanov numerical \cite{Rusanov1961} flux for one-dimensional problem is given as follows:
	\begin{equation}\label{3-num_rusanov}
		\textbf{F}_{rus} = \frac{1}{2}(\textbf{F}_L + \textbf{F}_R) - \frac{1}{2} S^+( \textbf{U}_R - \textbf{U}_L ),
	\end{equation}
	where $ S^+ $ is the positive wave speed given by Davis approximation \cite{Davis1988}:
	\begin{equation*}
		S^{+}=\max \left\{\left|u_{L}-c_{L}\right|,\left|u_{R}-c_{R}\right|,\left|u_{L}+c_{L}\right|,\left|u_{R}+c_{R}\right|\right\},
	\end{equation*}
	with the ``sound'' speed of the model:
	\begin{equation*}
		c^2 = gh + \frac{\lambda}{3}\left( \eta/h \right)^2.
	\end{equation*}
	Since contact characteristics are present in the system, we will also consider the HLLC Riemann solver proposed by Toro \cite{Toro2019}, \cite{Toro1994} and adapt it to the extended SGN model:
	\begin{equation*}
	\mathbf{F}_{hllc} = 
	\left\{
	\begin{array}{ll}
	{\mathbf{F}_L, } & {0 \leq S_{L}}, \\
	{\mathbf{F}_L^*,} & { S_{L} \leq 0 \leq S_{*}}, \\
	{\mathbf{F}_R^*, } & {S_{*} \leq 0 \leq S_{R}},\\
	{\mathbf{F}_R,} & { 0 \geq S_{R}}.
	\end{array}
	\right.
	\end{equation*}
	The intermediate fluxes are given by:
	\begin{equation*}
	\begin{aligned}
	& \textbf{F}_L^* = \textbf{F}_L + S_L(\textbf{U}_L^* - \textbf{U}_L),\\
	& \textbf{F}_R^* = \textbf{F}_R + S_R(\textbf{U}_R^* - \textbf{U}_R),
	\end{aligned}
	\end{equation*}
	where the intermediate conservative variables are:
	\begin{equation*}
	\textbf{U}_L^* = 
	\left(\begin{matrix}
	h_L^*\\
	h_L^*S^*\\
	h_L^*v_L\\
	h_L^*\eta_L\\
	h_L^*w_L
	\end{matrix}\right),
	\qquad
	\textbf{U}_R^* = 
	\left(\begin{matrix}
	h_R^*\\
	h_R^*S^*\\
	h_R^*v_R\\
	h_R^*\eta_R\\
	h_R^*w_R
	\end{matrix}\right),
	\end{equation*}
	with starred values defined by:
	\begin{equation*}
	h_L^* = h_L \frac{S_L - u_L}{S_L - S^*}, \qquad h_R^* = h_R \frac{S_R - u_R}{S_R - S^*}.
	\end{equation*}
	Here the middle wave speed is:
	\begin{equation*}
	S^* = \frac{p_R - p_L + h_R u_R (u_R - S_R) - h_L u_L(u_L-S_L)}{h_R(u_R-S_R) - h_L(u_L-S_L)},
	\end{equation*}
	\begin{equation*}
	S_{L}=\min \left\{u_{L}-c_{L}, u_{R}-c_{R}\right\}, \quad S_{R}=\max \left\{u_{L}+c_{L}, u_{R}+c_{R}\right\}.
	\end{equation*}
	
	The choice of the Riemann solver is motivated by the physical nature of the problem. If we consider smooth initial data like solitary waves, both Rusanov and HLLC fluxes produce the same results with no difference. However, when dealing with dispersive shock waves, HLLC keeps shock fronts sharper and better preserves the amplitudes of the trailing oscillations, thus we find it more suitable for the Riemann problem.

	The ODE part of \eqref{3-num_first-order-scheme} consists in resolution of the following subsystem:
	\begin{equation*}
		\pd{h}{t} = 0, \qquad \pd{u}{t} = 0, \qquad \pd{v}{t} = 0, \qquad \pd{\eta}{t} = hw, \qquad \pd{w}{t} = -\lambda\left(\frac{\eta}{h} - 1\right),
	\end{equation*}
	which admits the exact solution \cite{Favrie2017}:
	\begin{equation}\label{3-num_ode-exact-solution}
		\begin{aligned}
			& h_{ij}^{n+1} = h_{ij}^{n}, \qquad u_{ij}^{n+1} = u_{ij}^{n}, \qquad v_{ij}^{n+1} = v_{ij}^{n},\\
			& \eta_{ij}^{n+1} = h_{ij}^n + \left(\eta_{ij}^n - h_{ij}^n\right)\cos\left(\sqrt{\lambda}\frac{\Delta t}{h_{ij}^n}\right) + \frac{h_{ij}^n w_{ij}^n}{\sqrt{\lambda}} \sin\left(\sqrt{\lambda}\frac{\Delta t}{h_{ij}^n}\right),\\
			& w_{ij}^{n+1} = - \sqrt{\lambda}\Big( \frac{\eta_{ij}^n}{h_{ij}^n} - 1\Big)\sin\Big(\sqrt{\lambda}\frac{\Delta t}{h_{ij}^n}\Big) + w_{ij}^n \cos\Big(\sqrt{\lambda}\frac{\Delta t}{h_{ij}^n}\Big).
		\end{aligned}
	\end{equation}
	This exact solution defines the second step of \eqref{3-num_first-order-scheme}:
	\begin{equation*}
		\textbf{U}_{ij}^{(2)} = ODE_{exact} \left( \textbf{U}_{ij}^n \right) = h_{ij}^{n+1} \left(
			\begin{matrix}
			1\\
			u_{ij}^{n+1}\\
			v_{ij}^{n+1}\\
			\eta_{ij}^{n+1}\\
			w_{ij}^{n+1}
			\end{matrix}
		\right),
	\end{equation*}
	Eventually, we calculate the source terms from $ \textbf{U}_{ij}^{(2)} $ and utilize the explicit Euler procedure to update the numerical solution to the $ (n+1)^{\text{th}} $ layer:
	\begin{equation*}
		\textbf{U}_{ij}^{n+1} = \textbf{U}_{ij}^{(1)}  + \Delta t \, \textbf{S} \left(\textbf{U}_{ij}^{(2)} \right).
	\end{equation*}
	
\subsection{Second-order implicit-explicit method}

	The principal method we use for numerical resolution is the \textit{ARS(2,2,2)} implicit-explicit scheme \cite{Ascher1997}, \cite{Ascher1995}, \cite{Pareschi2001} of second order in space and time, which was already applied to one-dimensional hyperbolized dispersive systems in \cite{Dhaouadi2018}, \cite{Tkachenko2020a} and \cite{Richard2021}.
	\begin{equation}\label{3-num_ars222}
		\begin{aligned}
			& \textbf{U}_{ij}^{(1)} = \textbf{U}_{ij}^n + \alpha \Delta t \left( \mathcal{G}\big(\textbf{U}_{ij}^n\big) + \textbf{S}\big(\textbf{U}_{ij}^{(1)}\big) \right),\\
			& \textbf{U}_{ij}^{n+1} = \textbf{U}_{ij}^n + \Delta t \left( \delta\, \mathcal{G}\big(\textbf{U}_{ij}^n\big) + (1-\delta) \mathcal{G}\big(\textbf{U}_{ij}^{(1)}\big) \right) + \Delta t \left( \alpha \,\textbf{S}\big(\textbf{U}_{ij}^{n+1}\big) + (1-\alpha) \textbf{\textbf{S}}\big(\textbf{U}_{ij}^{(1)}\big) \right).
		\end{aligned}
	\end{equation}
	\begin{equation*}
		\alpha = 1 - \frac{1}{\sqrt{2}}, \qquad \delta = \alpha - 1.
	\end{equation*}
	The scheme consists of two steps, each of them containing two parts: the hyperbolic part which is solved explicitly and the implicit ODE part. Here $ \mathcal{G} $ is the same hyperbolic operator as \eqref{3-num_hyperbolic operator}, where the numerical fluxes are calculated using the MUSCL central difference piece-wise linear reconstruction, i.e. the left and right states of a one-directional Riemann problem are modified:
	\begin{equation*}
		\begin{aligned}
			& \textbf{U}_{L, i} = \textbf{U}_i - \frac{1}{2} \Delta_i,\\
			& \textbf{U}_{R, i} = \textbf{U}_i + \frac{1}{2} \Delta_i.
		\end{aligned}
	\end{equation*}
	The slope of an $ i $-th state is a pure central difference of the neighbor states taken without limiters:
	\begin{equation*}
		\Delta_i = \frac{1}{2}\left( \textbf{U}_{i+1} - \textbf{U}_{i-1}\right).
	\end{equation*}
	Then, the resolution of the Riemann problem is performed as for the first order method above with a Riemann solver of any choice:
	\begin{equation*}
		\textbf{F}_{i+\frac{1}{2}}^n = RP\left( \textbf{U}_{R, i}^n, \textbf{U}_{L, i+1}^n \right).
	\end{equation*}
	As we can notice, both implicit sub-steps are of the same form:
	\begin{equation*}
		\textbf{U}_{ij} = \textbf{U}_{ij}^0 + \alpha \,\Delta t \,\textbf{S}(\textbf{U}_{ij})
	\end{equation*}
	where $ \textbf{U}_{ij}^0 $ is known from the explicit calculations, and $ \textbf{U}_{ij} $ is an unknown to find. Luckily, this equation has an explicit solution:
	\begin{equation*}
		\begin{aligned}
			h_{ij} &= h_{ij}^0,\\
			u_{ij} &= u_{ij}^0,\\
			\eta_{ij} &= \frac{(h_{ij}^0)^2(\eta_{ij}^0+\alpha\Delta t w_{ij}^0) + \lambda \alpha^2 \Delta t^2 h_{ij}^0 }{(h_{ij}^0)^2 + \lambda \alpha^2 \Delta t^2},\\
			w_{ij} &= \frac{(h_{ij}^0)^2 w_{ij}^0 + \lambda \alpha \Delta t (h_{ij}^0 - \eta_{ij}^0)}{(h_{ij}^0)^2 + \lambda \alpha^2 \Delta t^2}.
		\end{aligned}
	\end{equation*}
	Hence, we use this solution for both sub-steps, taking $ \textbf{U}_{ij}^0 $ respectively for the first step:
	\begin{equation*}
	\textbf{U}_{ij}^0 = \textbf{U}_{ij}^n + \alpha \Delta t \, \mathcal{G}\big(\textbf{U}_{ij}^{n}\big),
	\end{equation*}
	and for the second one:
	\begin{equation*}
	\textbf{U}_{ij}^0 = \textbf{U}_{ij}^n + \Delta t \left( \delta\, \mathcal{G}\big(\textbf{U}_{ij}^n\big) + (1-\delta) \mathcal{G}\big(\textbf{U}_{ij}^{n+\frac{1}{2}}\big) \right) + (1-\alpha)\Delta t\, \textbf{\textbf{S}}\big(\textbf{U}_{ij}^{n+\frac{1}{2}}\big).
	\end{equation*}

	Stability studies of these numerical methods are non-trivial even for the first order case, and we use the standard 2-D CFL stability criteria relying on a common practice \cite{Godunov1976}:
	\begin{equation*}
	\frac{(|u|+c)\Delta t}{\Delta x} + \frac{(|v|+c)\Delta t}{\Delta y} \le CFL < 1.
	\end{equation*}
	Thus, the practical choice of the time step in numerical simulations is as follows:
	\begin{equation*}
	\Delta t=CFL / \max _{i, j}\left(\frac{\left|u_{i j}\right|+c_{i j}}{\Delta x}+\frac{\left|v_{i j}\right|+c_{i j}}{\Delta y}\right).
	\end{equation*}
	This criteria is slightly stronger than the directional maximum, as used by Colella \cite{Colella1990}. In practice, while testing dam break type problems with different CFL values, we noticed that it can reach slightly above $ 1 $ for the numerical solution not to explode.

%% file: sections/5_Numerical_results.tex
\subsection{One-dimensional}
	
	\subsubsection{Solitary wave}
\begin{figure}[H]
			\centering
			\includegraphics[width=\linewidth]{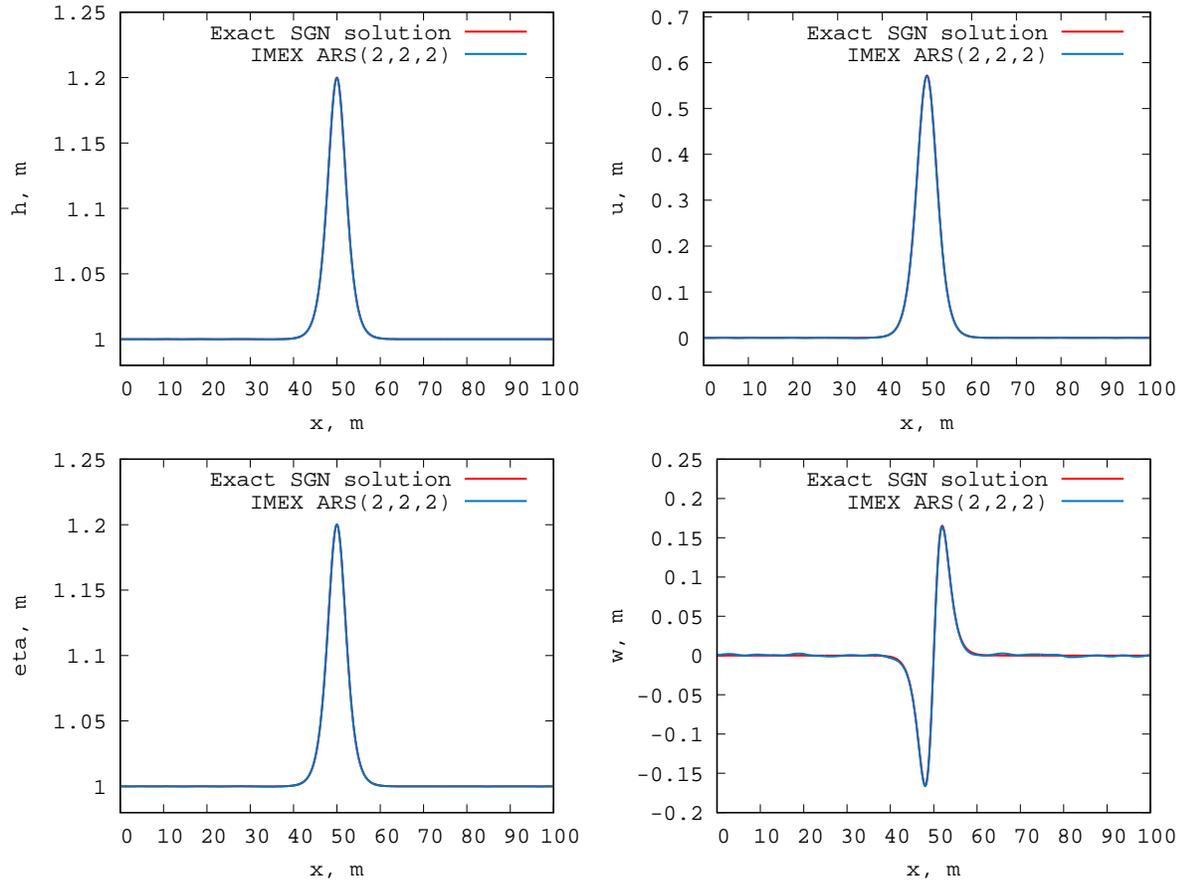}
			\caption{Numerical and exact solitary wave solutions at $ T = 100/D = 29.14573 \, s $. Rusanov flux, $ CFL=0.9 $.}
			\label{fig:4_res_oned-dumbser-ord2}
		\end{figure}

		Authors in \cite{Favrie2017} considered solitary waves as primary validation tests for the first order splitting. We will complement these studies with classical 1-D solitary wave tests using the second-order IMEX method. For the two following problems we will use the Rusanov solver to calculate the numerical flux in \eqref{3-num_ars222}. We now consider the propagation of the original Serre-Green-Naghdi solitary wave given by:
		\begin{equation*}
		\begin{aligned}
		& h(x, t) = h_0 + a\,\mathrm{sech}^2 \big(\kappa (x-Dt)\big),\\
		& u(x,t) = D\left( 1 - \frac{h_0}{h(x,t)} \right),
		\end{aligned}
		\end{equation*}
		where
		\begin{equation*}
		\kappa = \sqrt{\frac{3a}{4h_0^2 (h_0 + a)}}, \qquad D = \sqrt{g(h_0+a)}.
		\end{equation*}
		We take $ h = \eta $ at $ t = 0 $ to satisfy the initial equilibrium condition. In addition, as remarked in \cite{Duchene2019}, it is important to couple the initial data for $ w $ with its definition, i.e. $ w = \dot{\eta} $ and thus $ w = \dot{h} $:
		\begin{equation*}
		\left.w\right|_{t=0} = \left.(h_t + u h_x)\right|_{t=0} = \left.-hu_x\right|_{t=0} =  - D \frac{h_0}{h(x,0)} h_x(x,0).
		\end{equation*}
		Here we take $ h_0 = 1\,m $, $ a = 0.2\,m $, $ g = 9.81\,m/s^2 $ and $ \lambda = 1200\, m^2/s^2 $. The domain is $ 100 $ meters long, we take 2000 mesh points, impose periodic boundary conditions and run the calculations until the full period is reached, i.e. the final simulation time is $ T = 100/D $. The results obtained with the IMEX scheme are shown in the Fig. \ref{fig:4_res_oned-dumbser-ord2}. The numerical solution is very close to the exact one except for some small-amplitude tailing oscillations of the variable $ w $.
		
	\subsubsection{Soliton head-on collision}

\begin{figure}[H]
			\centering
			\includegraphics[width=\linewidth]{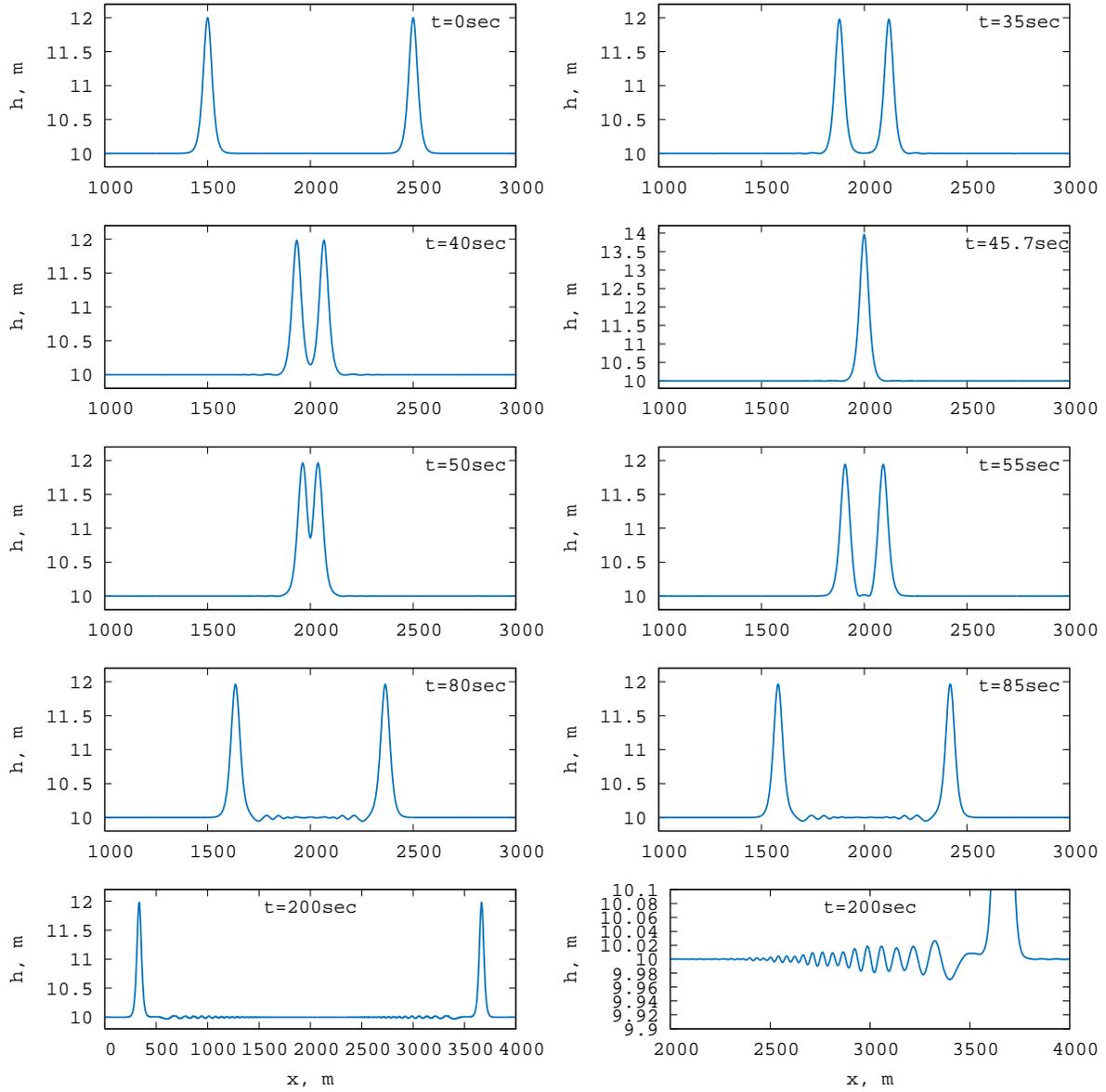}
			\caption{Solitary wave head-on collision at $ t = 200\,s $. Rusanov flux, $ CFL=0.9 $.}
			\label{fig:4_res_oned-collision-moments}
		\end{figure}

\begin{figure}[H]
			\centering
			\includegraphics[width=\linewidth]{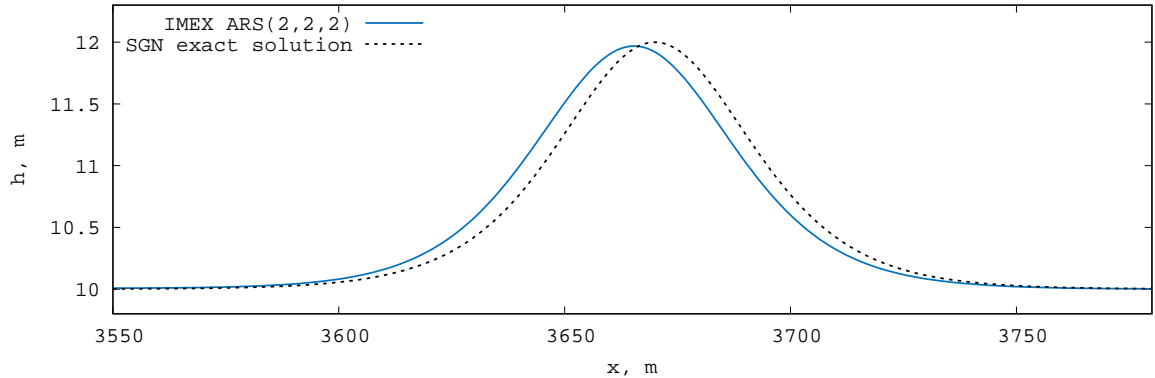}
			\caption{Comparison of a single exact SGN solitary wave to the numerical solution to the head-on collision problem.}
			\label{fig:4_res_oned-collision}
		\end{figure}

		We reproduce another classic test, notably the head-on collision of two solitary waves of equal amplitude described in \cite{Mitsotakis2014} and \cite{Duran2015}, where the authors used the finite element Galerkin/finite element discretization for different versions of Serre-Green-Naghdi equations. Initially two solitary waves are placed at $ x = 1500\,m $ and $ x = 2500\,m $ in a $ 4000 \, m$ long domain and directed towards each other. We take $ h_0 = 10\,m $, $ a = 2\,m $, $ g = 9.81\,m/s^2 $, $ \lambda = 2400\, m^2/s^2 $, $ 4000 $ mesh points, and the final time $ T = 200s $. The solution is pictured on several snapshots in the Fig. \ref{fig:4_res_oned-collision-moments}. Small amplitude oscillations follow the solitary waves which is clearly seen on the last zoomed-in section corresponding $ t = 200\,s $. One can also observe a small phase shift and slight amplitude loss as compared to a single traveling solitary wave in the same setup (Fig. \ref{fig:4_res_oned-collision}). Thus, the finite volume and finite element methods give out the same results as expected.
	
	\subsubsection{Dam break problem}
\begin{figure}[H]
			\centering
			\includegraphics[width=\linewidth]{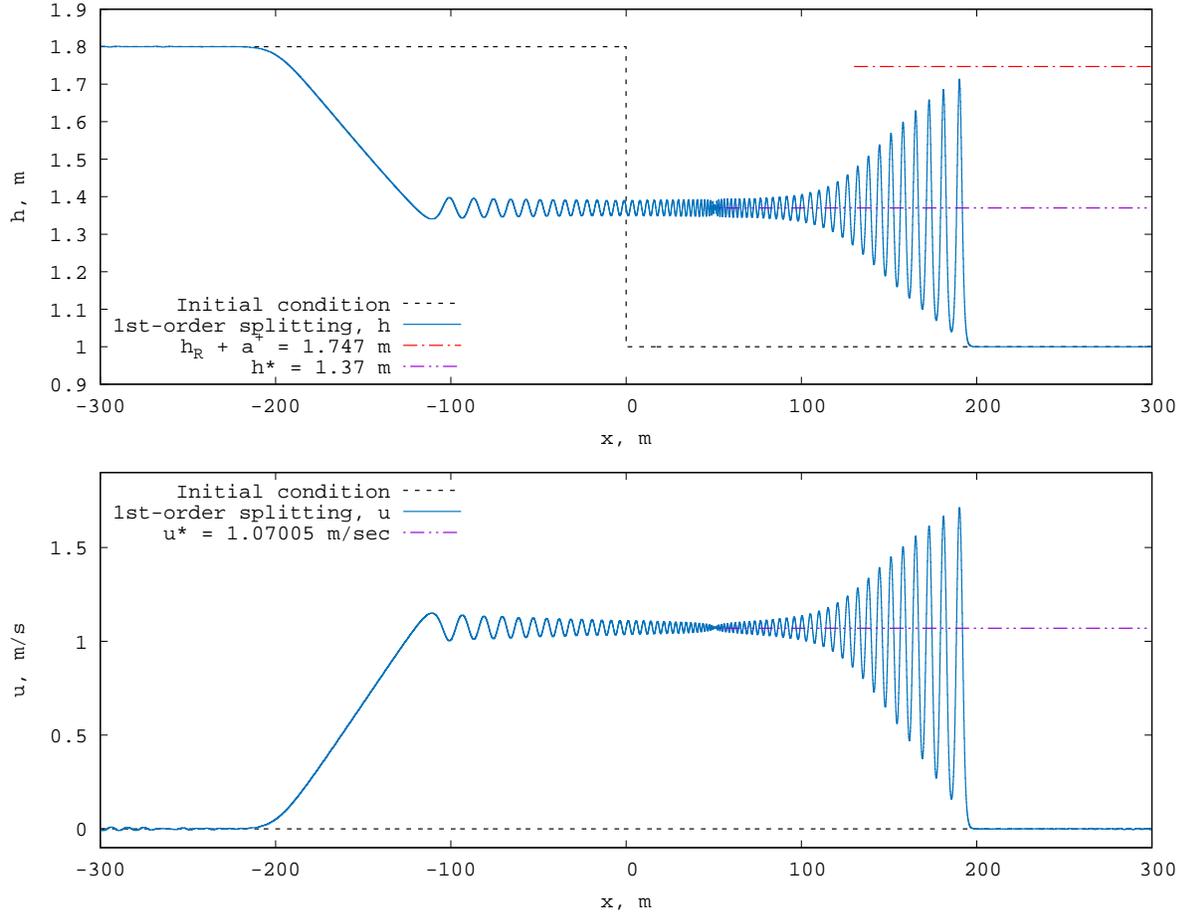}
			\caption{Numerical solution to 1-D dam break problem \ref{4-res_ic_dam-break-1d} using the first-order splitting: water depth and horizontal velocity profiles at $ t = 47.434\,s $. HLLC flux, 768000 mesh points, $ CFL = 0.8 $.}
			\label{fig:4_res_dam-break-1d_1st-order}
		\end{figure}

		Let us consider the propagation of dispersive shock waves to demonstrate the robustness of the second-order method. For all the following tests we will use the HLLC Riemann solver to calculate the numerical flux, as explained in the previous chapter. Initially two different states are separated by an infinitely thin barrier. This configuration is imitated by a piece-wise constant initial data with a discontinuity at $ x = 0 $:
		\begin{equation}\label{4-res_ic_dam-break-1d}
		h(x, 0) = \left\{
		\begin{aligned}
		& h_L, & x \le 0, \\
		& h_R, & x > 0, \\
		\end{aligned}
		\right. \qquad u(x,0) = 0.
		\end{equation}
		As we did in the previous example, we add the corresponding initial data for $ \eta $ and $ w $:
		\begin{equation*}
		\eta(x, 0) = h(x, 0), \qquad w(x, 0) = \left.-hu_x\right|_{t=0} = 0.
		\end{equation*}
		The second-to-last equality is a corollary of the 1-D mass conservation law.  Initially we consider the problem using the first-order splitting. We take $ 768000 $ mesh points, $ g = 9.81\,m/s^2 $, $ \lambda = 300\,m^2/s^2 $, $ h_L = 1.8\,m $ and $ h_R = 1\,m $. At $ t = 0 $ the barrier is removed and the initial configuration breaks up: the discontinuity divides into a dispersive shock wave propagating to the right and a rarefaction wave to the left, leaving the plateau region between them, see Fig. \ref{fig:4_res_dam-break-1d_1st-order}.
	
		The form of the solution corresponding the discontinuous initial data, as well as the amplitude of the leading solitary wave, is in a good agreement with the results of the same test for the original dispersive SGN model, obtained in  \cite{Gavrilyuk2020} using a semi-discrete finite method \cite{Ketcheson2008}, \cite{Ketcheson2013}, with the smoothed initial data:
		\begin{equation*}
		h(x, 0) = h_R + \frac{h_L-h_R}{2}\left(1-\tanh\left(\frac{x}{\alpha}\right)\right), \qquad u(x,0) = 0,
		\end{equation*}
		and $ \alpha = 0.4 $ (i.e. the measure of the transition region between the left and the right states could be considered negligible). The analysis of Riemann invariants of the shallow-water system, coupled with the analysis of Witham system for Serre-Green-Naghdi equations \cite{El2006}, \cite{Whitham1974}, namely its solutions of simple wave type, allow to recover the approximate values $ (h^*,u^*) $ of the mean flow dividing the rarefaction wave and the dispersive shock zones:
		\begin{equation}\label{4-res_el-h-_u-}
		h^* = \frac{(\sqrt{h_L}+\sqrt{h_R})^2}{4}, \qquad u^* = 2\big(\sqrt{g h^*}-\sqrt{g h_R}\big).
		\end{equation}
		We can see that those quantities match the numerical values of these parameters (purple dashed double dotted lines in Fig. \ref{fig:4_res_dam-break-1d_1st-order}. The second-order asymptotic approximation of the amplitude of the lead soliton $ a^+ $ \cite{El2005} is:
		\begin{equation}\label{4-res_el-amplitude}
			a^+ = \delta_0 - \frac{1}{12} \delta_0^2 + O(\delta_0^3),
		\end{equation}
		where $ \delta_0 $ denotes the initial jump value. The numerical value of $ a^+ $ is also in a good agreement with the approximate expression (red dashed single dotted line in Fig.  \ref{fig:4_res_dam-break-1d_1st-order}). Although the results are rather accurate, the major disadvantage is still the large computational time due to large number of cells. We developed an MPI parallel algorithm for the first-order method to perform such simulations: we used $ 48 $ $ 2.3 $ GHz processors for this test, which took 3 hours and 15 minutes of calculations.
		
\begin{figure}[H]
			\centering
			\includegraphics[width=\linewidth]{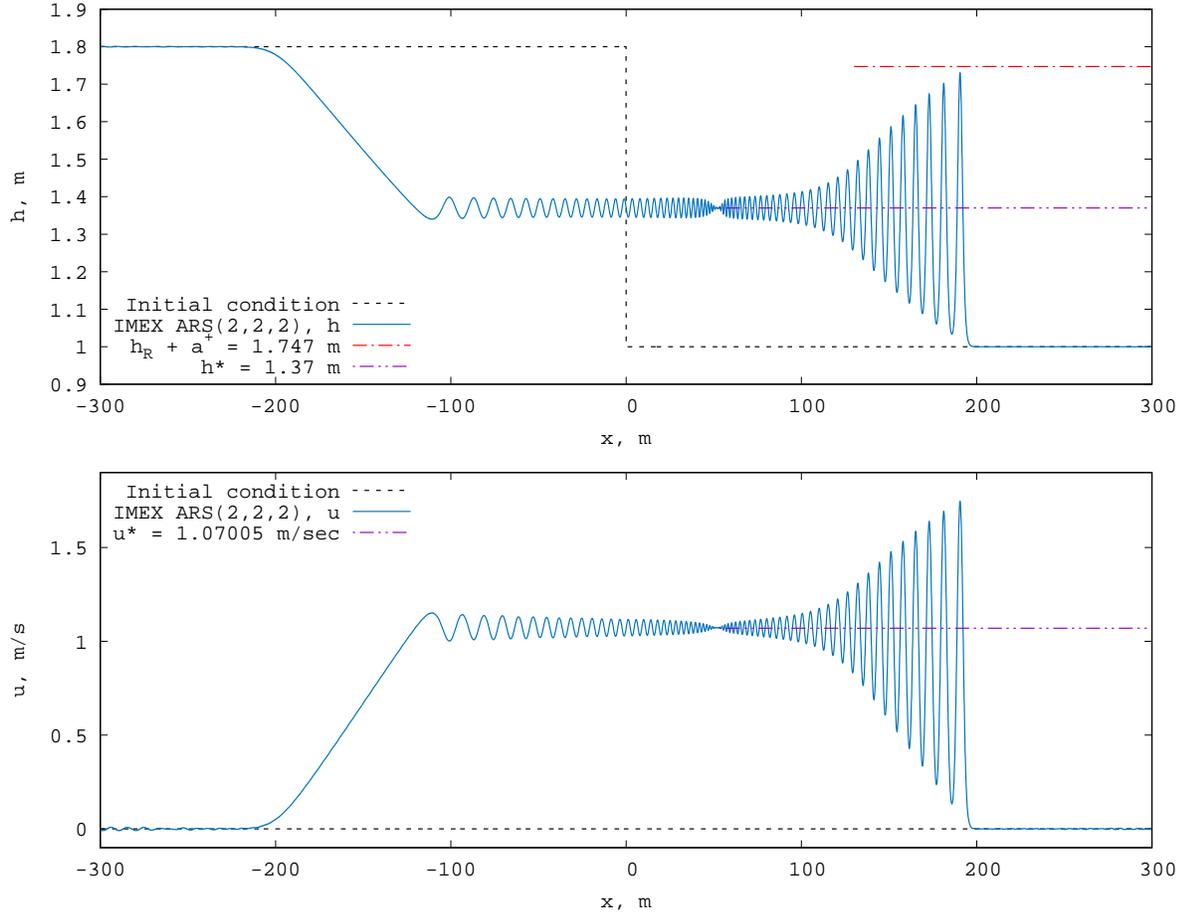}
			\caption{Numerical solution to 1-D dam break problem \ref{4-res_ic_dam-break-1d} using the IMEX method: water depth and horizontal velocity profiles at $ t = 47.434\,s $. HLLC flux, $ 8000 $ mesh points, $ CFL = 0.95 $.}
			\label{fig:4_res_dam-break-1d_2nd-order}
		\end{figure}

		Now consider the numerical solution to the same problem using the second-order IMEX method taking $ 8000 $ mesh points (see Fig. \ref{fig:4_res_dam-break-1d_2nd-order}).
		One can notice that asymptotic parameters are in better agreement with \eqref{4-res_el-h-_u-} and \eqref{4-res_el-amplitude} than those obtained with a first-order method, and 96 times less cells were needed. Since the number of points significantly decreased, we did not parallelize the code and used only one 2.3 GHz processor, which took 29 seconds of calculations to reproduce an even more accurate result. The calculation time and the processor information used in tests for both methods are resumed in the Table \ref{tab:4_res_cpu_dam-break-1d}.
		\begin{table}[H]
			\resizebox{\textwidth}{!}{%
				\begin{tabular}{|c|c|c|c|c|}\hline
					Method          & Processor & Number of processors & Mesh points  & Calculation time                                 \\ \hline
					First-order splitting                & Intel (R) Xeon(R) CPU E7-4850 v2 @ 2.30 GHz  & 48  & 768000     & 3h 15m 41s\\ \hline
					IMEX ARS(2,2,2)                   & Intel (R) Core (TM) i5-7360U CPU @ 2.30 GHz  & 1     & 8000     & 29s   \\ \hline
				\end{tabular} %
			}
			\caption{Calculation time and processor information.\label{tab:4_res_cpu_dam-break-1d}}
		\end{table}
		Thus, the IMEX method demonstrates better precision and demands much less computation resources to reproduce a dispersive shock wave.
	
\subsection{Two-dimensional}
	\subsubsection{A symmetrical 2-D dam break problem} 
\begin{figure}[H]
			\centering
			\begin{subfigure}{.5\textwidth}
				\centering
				\includegraphics[width=\textwidth]{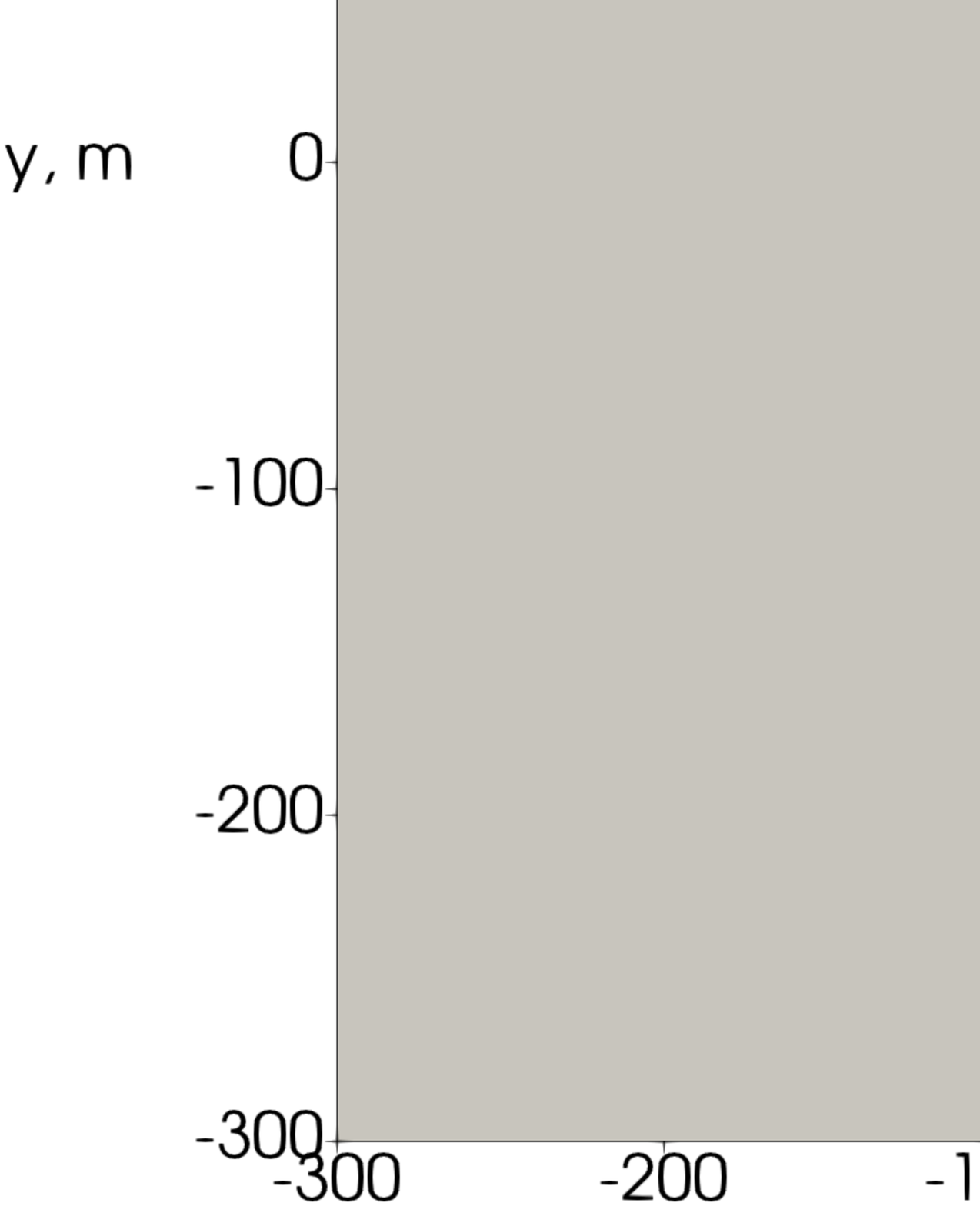}
			\end{subfigure}%
			\begin{subfigure}{.5\textwidth}
				\centering
				\includegraphics[width=\textwidth]{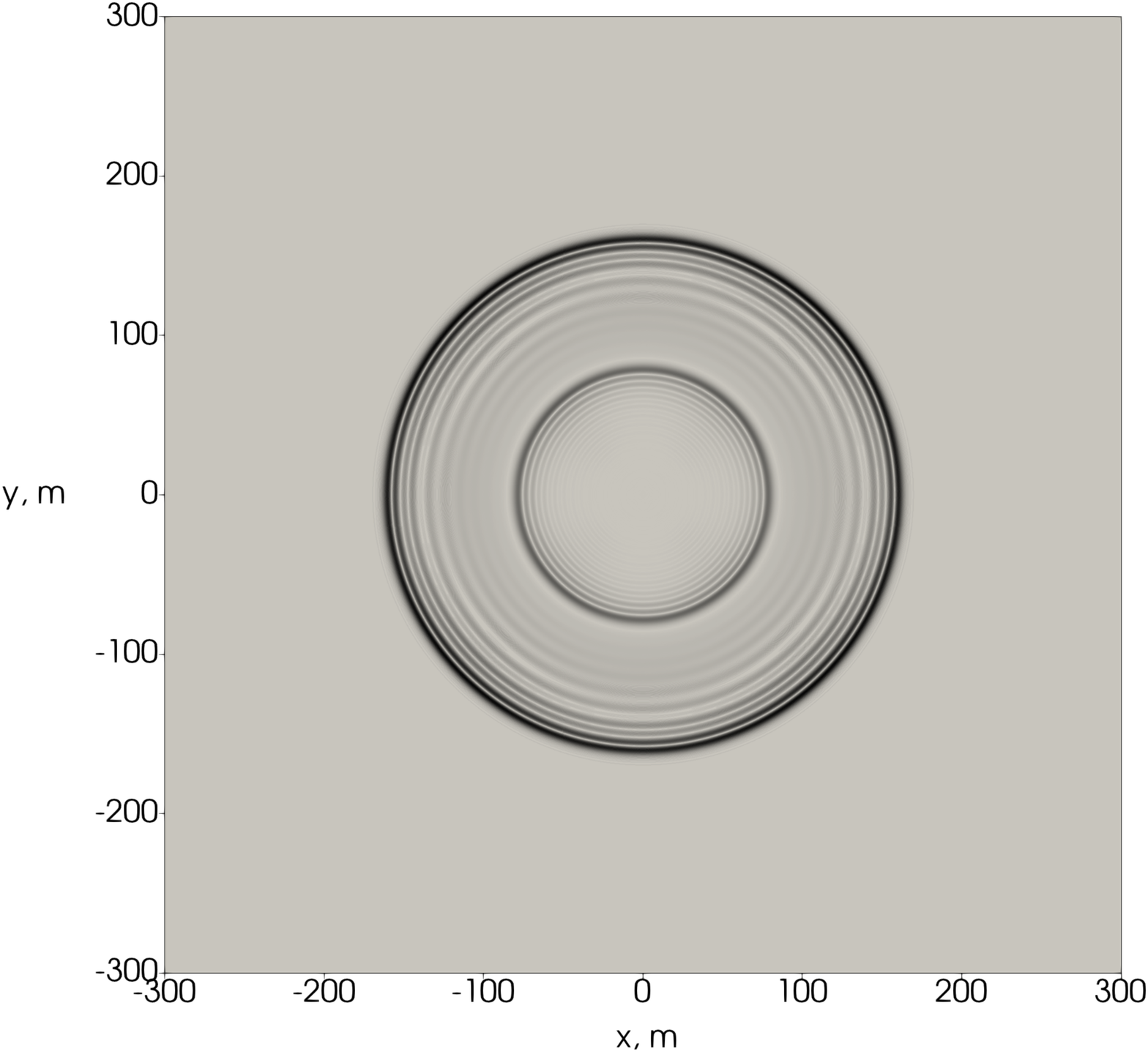}
			\end{subfigure}
			\caption{Numerical Schlieren visualization of the initial and final contour plots of the water depth corresponding to the 2-D dam break problem \eqref{4-res_ic_dam-break-2d-circle} at $ t=40\,s $, using the first-order splitting. HLLC flux, $ CFL = 0.5 $. The Schlieren function is taken $ \ln\left( 1 + 2 \abs{\nabla h} \right) $.}
			\label{fig:4_res_rp-2d-godunov-circle}
		\end{figure}

\begin{figure}[H]
			\centering
			\includegraphics[width=\textwidth]{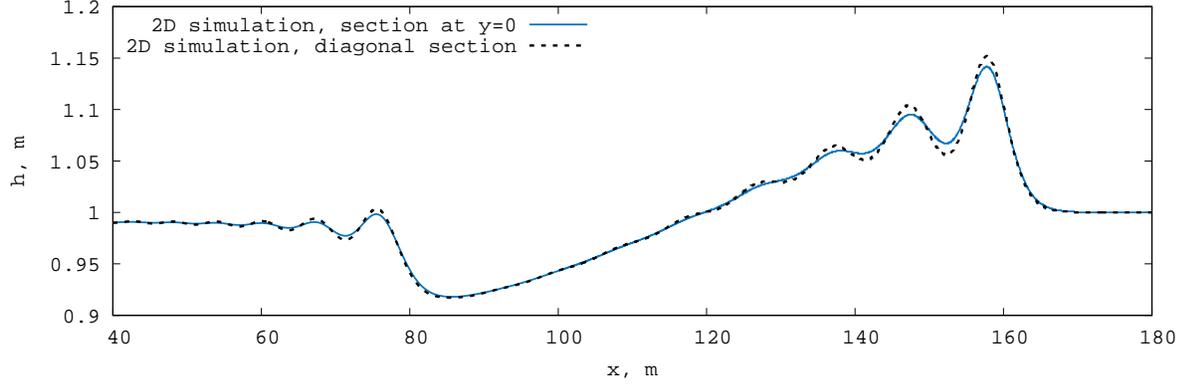}
			\caption{Two different cross-sections of Fig. \ref{fig:4_res_rp-2d-godunov-circle}: at the line $ y = 0 $ and the diagonal one.}
			\label{fig:4_res_rp-2d-godunov-circle_section0_vs_section45}
		\end{figure}

\begin{figure}[H]
			\centering
			\includegraphics[width=\textwidth]{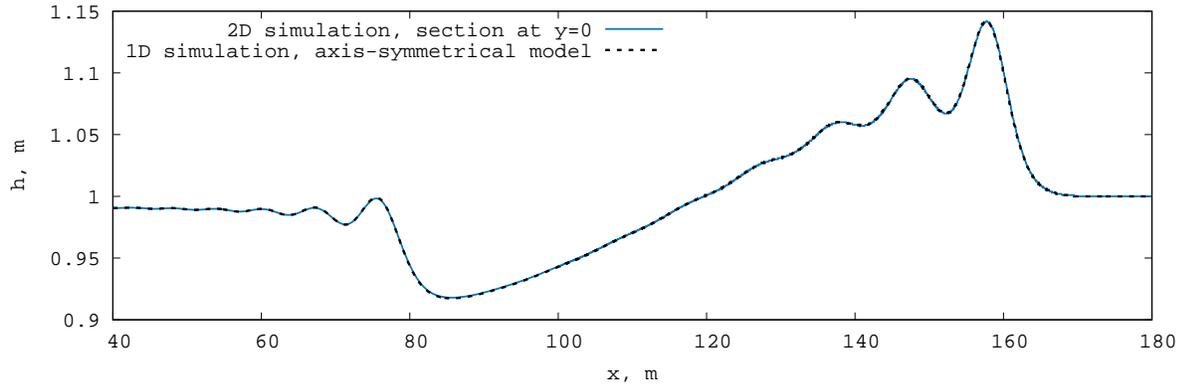}
			\caption{Comparison of the 2-D numerical solution to \eqref{4-res_ic_dam-break-2d-circle} (cross-section at $ y = 0 $) to the one-dimensional one obtained with the axis-symmetrical analogue \eqref{4-res_ext-sgn-polar-1d-conservative_2}.}
			\label{fig:4_res_rp-2d-godunov-circle_polar-1d_vs_section0}
		\end{figure}

		Consider a two-dimensional dam break problem: we impose a piecewise-constant initial data as follows: a circle of radius $ R_c $ is placed in the center of the computational domain $ {-300~m < x,y < 300~m} $, the water depth is $ h_{in} = 1.8\,m $ inside the circle and $ h_{out} = 1.0\,m $ outside, and the initial velocity is zero:
		\begin{equation}\label{4-res_ic_dam-break-2d-circle}
			h(x, 0) = \eta(x, 0) = \left\{
			\begin{aligned}
				& h_{in}, & x^2 + y^2 \le R_c^2, \\
				& h_{out}, & x^2 + y^2 > R_c^2. \\
			\end{aligned}
			\right., \qquad u(x,y,0) = v(x,y,0) = w(x,y,0) = 0,
		\end{equation}
		with $ g = 9.81\,m/s^2 $, $ \lambda = 75\,m^2/s^2 $ and $ 10000 \times 10000 $ mesh points. The numerical solution at $ t = 40\,s $ obtained with the first-order method \eqref{3-num_first-order-scheme} and HLLC flux is shown in the Fig. \ref{fig:4_res_rp-2d-godunov-circle}.
		Although the problem is symmetric, the cross-sections taken at different axes, namely at $ y = 0 $ and at the diagonal axis from $ (-300~m, -300~m) $ to $ (300~m, 300~m) $, slightly differ from each other, since we use the Cartesian mesh, which is not adapted to the symmetrical nature of the original problem (see Fig. \ref{fig:4_res_rp-2d-godunov-circle_section0_vs_section45}).
		In order to single out the correct solution, we perform the same test for a 1-D axis-symmetric version of \eqref{3-num_ext-sgn-2d-conservative-vector} in polar coordinates:
		
		\begin{equation}\label{4-res_ext-sgn-polar-1d-conservative_2}
		\begin{aligned}
		& \pd{h}{t} + \pd{(hu)}{r} = -\frac{hu}{r}, \qquad r = \sqrt{x^2+y^2},\\
		& \pd{(hu)}{t} + \pd{(hu + p)}{r} = -\frac{hu}{r},\\
		& \pd{(h\eta)}{t} + \pd{(hu\eta)}{r} = hw - \frac{hu\eta}{r},\\
		& \pd{(hw)}{t} + \pd{(huw)}{r} = -\lambda\left(\frac{\eta}{h}-1\right) - \frac{huw}{r}.
		\end{aligned}
		\end{equation}
		The numerical solution to the 1-D version of the considered dam break problem perfectly corresponds the cross-section at axis $ y = 0 $ of the full 2-D solution to \eqref{4-res_ic_dam-break-2d-circle} which is thus the correct one (see Fig. \ref{fig:4_res_rp-2d-godunov-circle_polar-1d_vs_section0}).

	\subsubsection{A non-radial 2-D dam break problem}
\begin{figure}[H]
			\centering
			\begin{subfigure}{.5\textwidth}
				\centering
				\includegraphics[width=\textwidth]{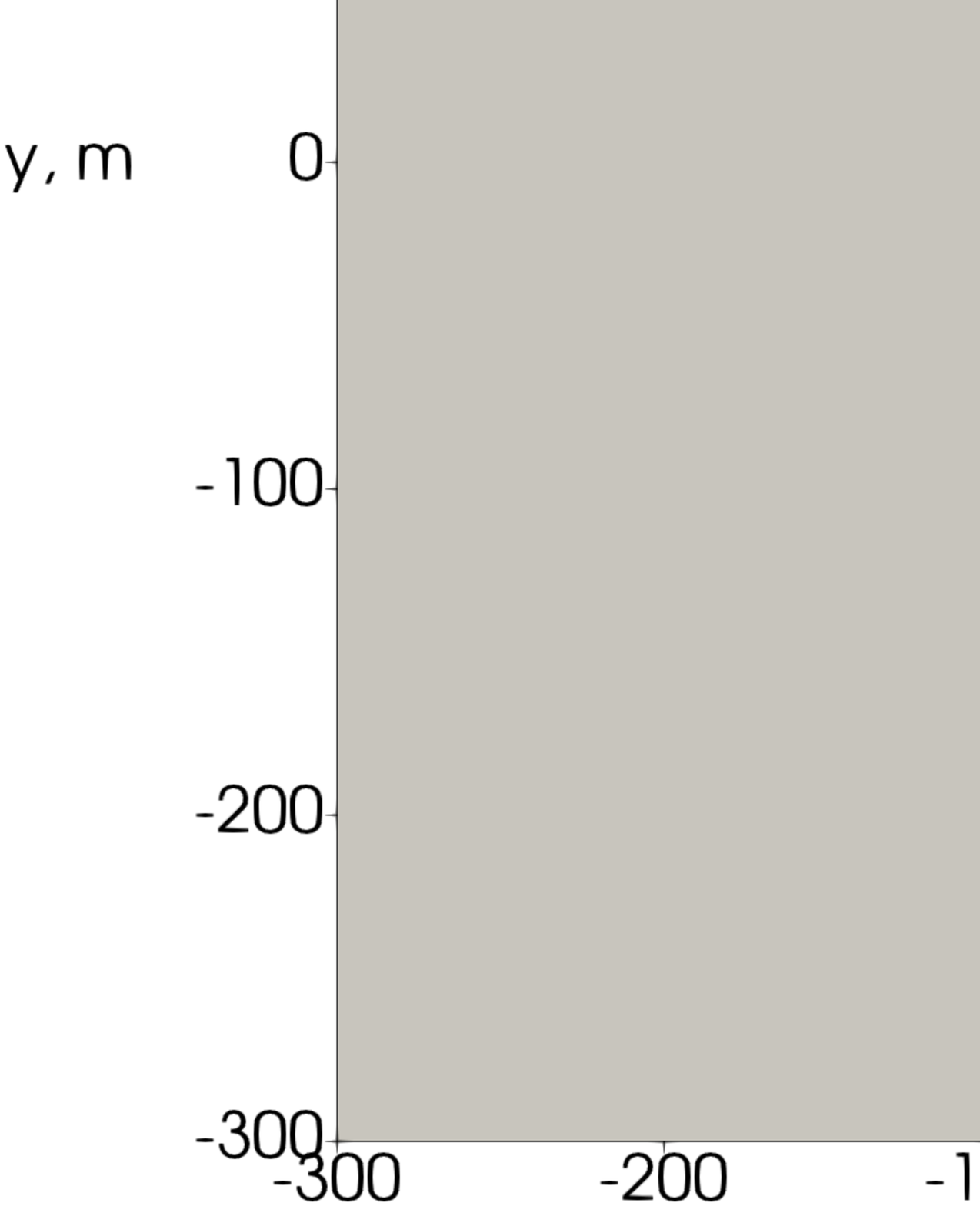}
			\end{subfigure}%
			\begin{subfigure}{.5\textwidth}
				\centering
				\includegraphics[width=\textwidth]{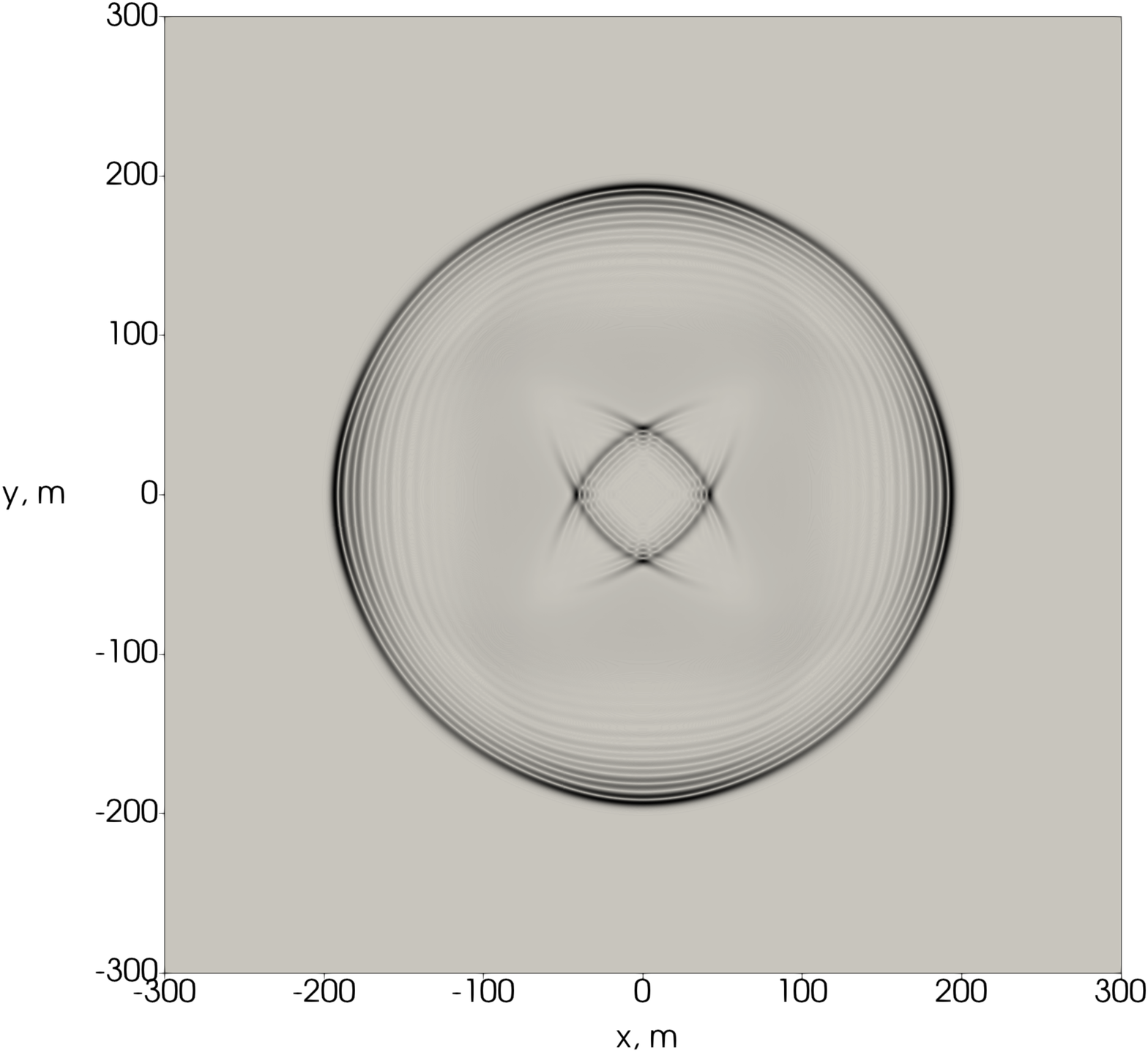}
			\end{subfigure}
			\caption{Initial condition and numerical solution to the non-symmetrical 2-D dam break problem \eqref{4-res_ic_dam-break-2d-square} at $ t = 40\,s $ using the first-order splitting. HLLC flux, $ CFL = 0.5 $.}
			\label{fig:4_res_rp-2d-godunov-square}
		\end{figure}

\begin{figure}[H]
			\centering
			\includegraphics[width=\textwidth]{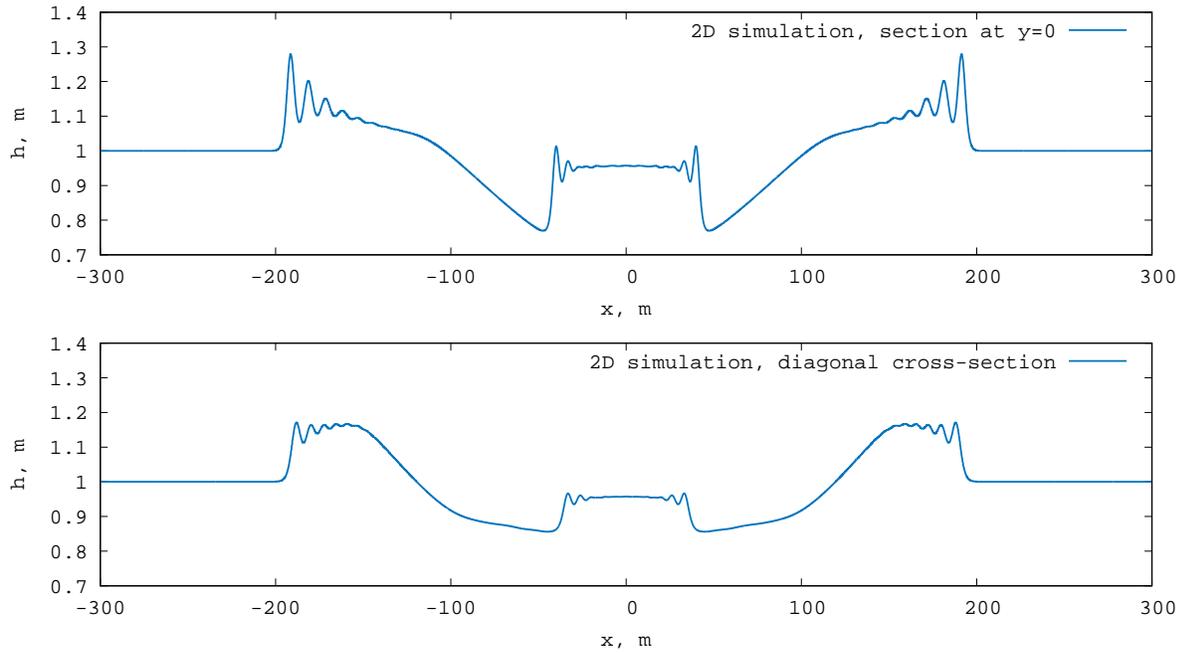}
			\caption{Horizontal and diagonal cross-sections of Fig. \ref{fig:4_res_rp-2d-godunov-square}.}
			\label{fig:4_res_rp-2d-godunov-square_section}
		\end{figure}

		In addition to the previous one, consider a non-radial Riemann problem. The initial configuration is the same as in \eqref{4-res_ic_dam-break-2d-circle} but the elevated surface is of square form (see Fig. \ref{fig:4_res_rp-2d-godunov-square}, left):
		\begin{equation}\label{4-res_ic_dam-break-2d-square}
			h(x, y, 0) = \eta(x, y, 0) = \left\{
			\begin{aligned}
				& h_{in}, & |x| \le d_s/2~\mathrm{and}~|y| \le d_s/2,\\
				& h_{out}, & \mathrm{otherwise}. \\
			\end{aligned}
			\right., \qquad	u(x,y,0) = v(x,y,0) = w(x,y,0) = 0.
		\end{equation}
		where $ d_s $ is a square side length. The water depth is $ h_{in} = 1.8~m $ inside the square and $ h_{out} = 1.0~m $ outside. The square side is $ d_s = 80~m $, $ g = 9.81\,m/s^2 $ and $ \lambda = 75\,m^2/s^2 $.	One can observe the non-symmetrical structures qualitatively different from the symmetrical case in the Fig. \ref{fig:4_res_rp-2d-godunov-square_section}.
		
	\subsubsection{2-D dam break problems with the second order method}
\begin{figure}[H]
			\centering
			\begin{subfigure}{.5\textwidth}
				\centering
				\includegraphics[width=\textwidth]{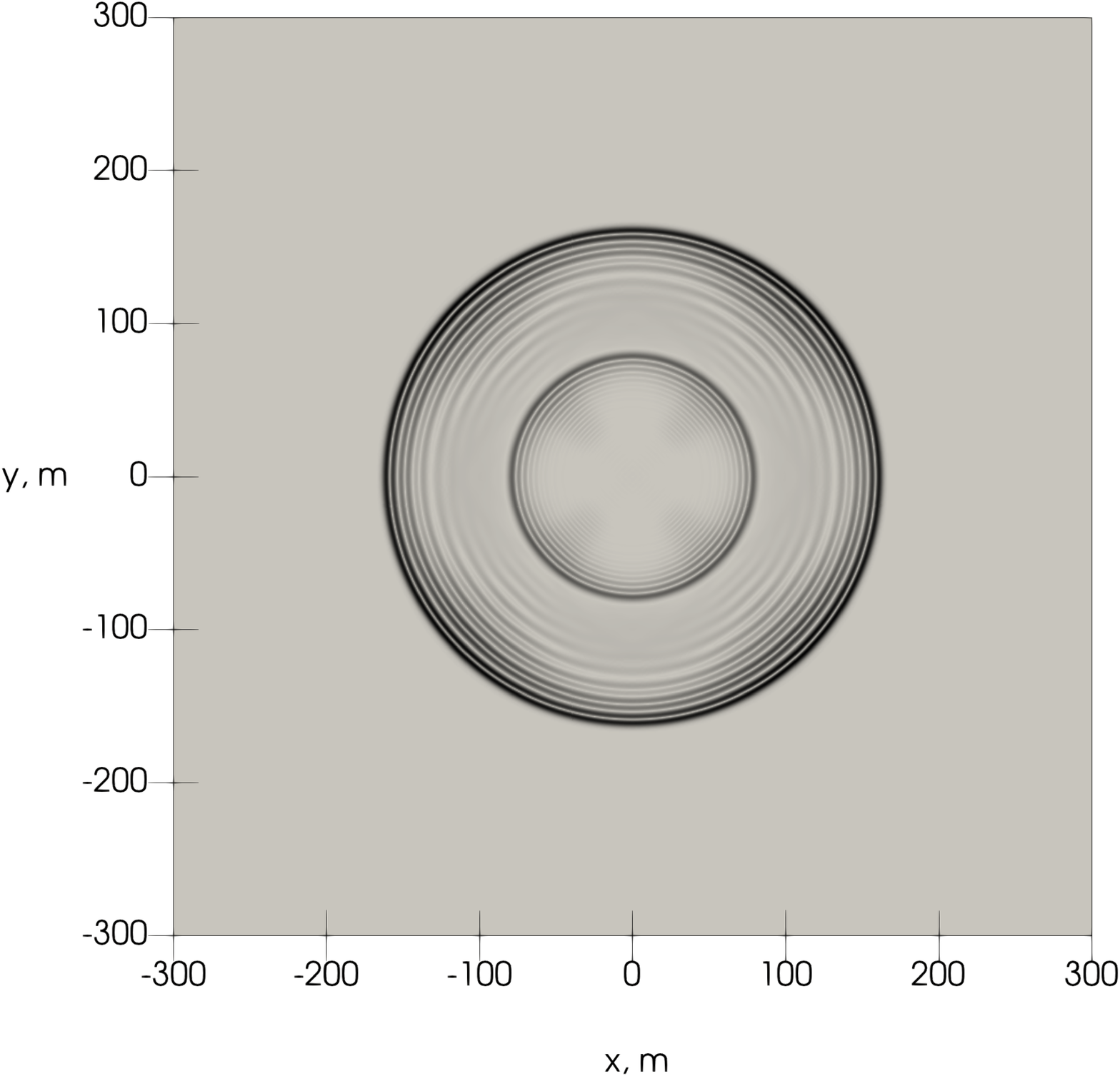}
				\caption{Circular cylinder}
				\label{fig:4_res_rp-2d-imex-circle}
			\end{subfigure}%
			\begin{subfigure}{.5\textwidth}
				\centering
				\includegraphics[width=\textwidth]{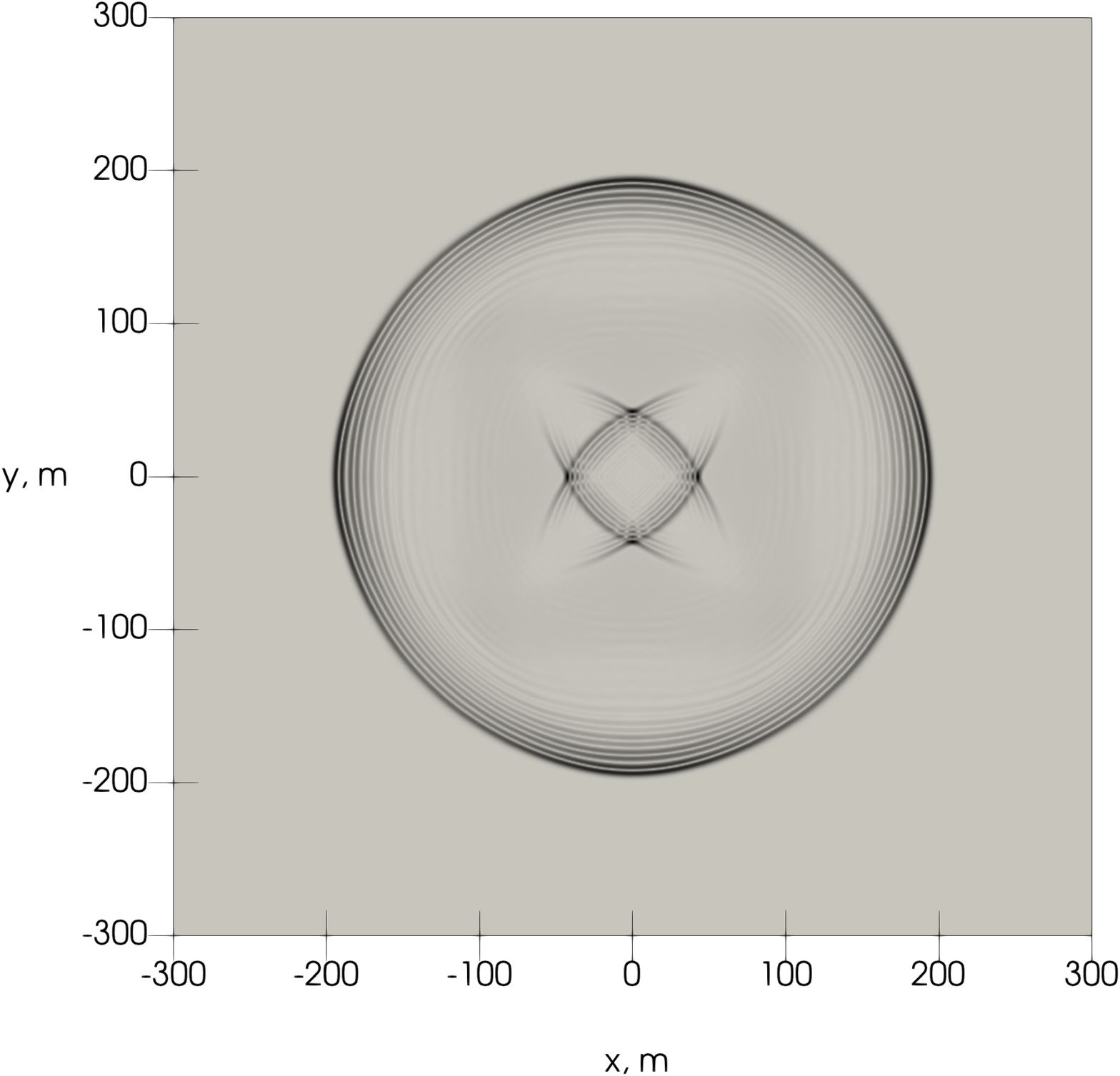}
				\caption{Square cylinder}
				\label{fig:4_res_rp-2d-imex-square}
			\end{subfigure}
			\caption{Numerical solution to the 2-D dam break problems \eqref{4-res_ic_dam-break-2d-circle} and \eqref{4-res_ic_dam-break-2d-square} at $ t = 40\,s $ using the IMEX method. HLLC flux, $  800\times800 $ mesh points, $ CFL = 0.5 $.}
			\label{fig:4_res_rp-2d-imex}
		\end{figure}

\begin{figure}[H]
			\centering
			\includegraphics[width=\textwidth]{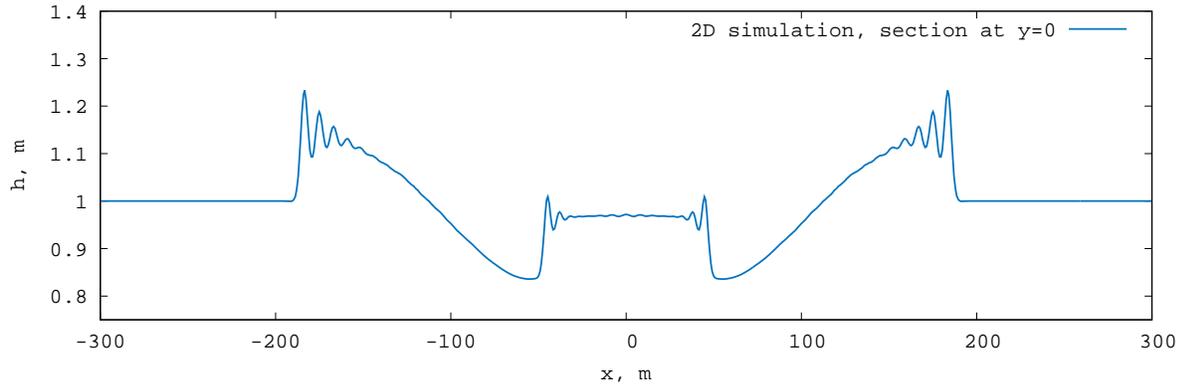}
			\caption{Horizontal cross-section of Fig. \ref{fig:4_res_rp-2d-imex-circle}.}
			\label{fig:4_res_rp-2d-imex-circle-section}
		\end{figure}

\begin{figure}[H]
			\centering
			\includegraphics[width=\textwidth]{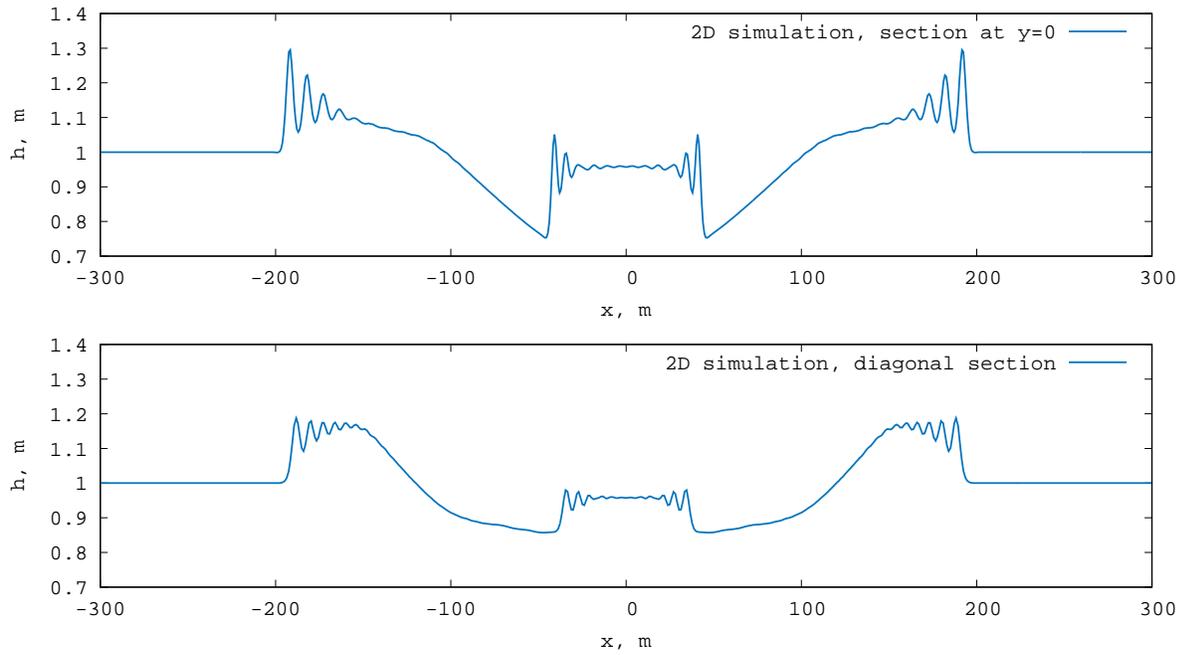}
			\caption{Horizontal and diagonal cross-sections of Fig. \ref{fig:4_res_rp-2d-imex-square}.}
			\label{fig:4_res_rp-2d-imex-square-sections}
		\end{figure}

		Let us consider the same 2D problems \eqref{4-res_ic_dam-break-2d-circle}, \eqref{4-res_ic_dam-break-2d-square} and use the second-order IMEX method \eqref{3-num_ars222} with HLLC solver, trying to achieve the results similar to those obtained with the first-order splitting. It turns out that only $ 800 \times 800 $ points is needed which is around 156 times less than we used in the previous example. The results are presented in Figs. \ref{fig:4_res_rp-2d-imex} -- \ref{fig:4_res_rp-2d-imex-square-sections}. The calculation time and the processor information used in both tests \eqref{4-res_ic_dam-break-2d-circle} and \eqref{4-res_ic_dam-break-2d-square} are summarized in the Table \ref{tab:4_res_cpu_dam-break-2d}.
		
		\begin{table}[H]
			\resizebox{\textwidth}{!}{%
				\begin{tabular}{|c|c|c|c|c|}\hline
					Method &Test & Number of processors & Calculation time & Processor \\ \hline
					First-order splitting & Circle 10000x10000 & 48 & 5h 12m 58s & Intel (R) Xeon(R) CPU E7-4850 v2 @ 2.30 GHz  \\ \hline
					First-order splitting & Square 10000x10000 & 48 & 5h 39m 7s  & Intel (R) Xeon(R) CPU E7-4850 v2 @ 2.30 GHz  \\ \hline
					IMEX ARS(2,2,2) & Circle 800x800 & 1 & 16m 25s & Intel (R) Core (TM) i5-7360U CPU @ 2.30 GHz  \\ \hline
					IMEX ARS(2,2,2) & Square 800x800 & 1 & 17m 28s & Intel (R) Core (TM) i5-7360U CPU @ 2.30 GHz  \\ \hline
					
				\end{tabular} %
			}
			\caption{Calculation time and processor information.\label{tab:4_res_cpu_dam-break-2d}}
		\end{table}

%% file: sections/6_Conclusion.tex
We studied multidimensional nonlinear dispersive models describing, in particular, shallow water flows and bubbly fluids. They are Euler-Lagrange equations for a Lagrangian depending on state variables and their first material derivatives. Using the extended Lagrangian approach proposed in \cite{Favrie2017}, we derived a Galilean invariant and unconditionally hyperbolic system which approximates the corresponding physical models.

To perform numerical simulations we consider the multi-D SGN equations. The robust $ ARS(2, 2, 2) $ IMEX method was used. It requires little mesh points to reach a good precision of the numerical solutions. The numerical results are in good agreement with the available exact solutions and those obtained with other numerical methods.

\subsection*{Acknowledgments}
Authors would like to thank Frederic Couderc, Firas Dhaoudi and Nicolas Favrie for helpful suggestions and discussions.

%% file: sections/Appendix.tex
\section{Extended model derivation and study}\label{appendix:extended-model}
	\subsection{Hamilton's principle}\label{appendix:extended-model-hamilton}
		Consider the following Lagrangian:
		\begin{equation*}
			\mathcal{L} = \int\limits_{D(t)} \left(\frac{ \rho\left| \textbf{u} \right|^2 }{ 2 } - W(\rho,\eta,\dot{\eta}) \right)d D.
		\end{equation*}
		The corresponding action functional is
		\begin{equation*}
			a = \int\limits_{t_0}^{t_1} \mathcal{L} \,dt = \int\limits_{t_0}^{t_1} \int\limits_{D(t)} \left( \frac{\rho\left|\textbf{u}\right|^2}{2} - W(\rho,\eta,\dot{\eta})\right)dD\,dt.
		\end{equation*}
		The state of the system is characterized by two variables $ \rho $ and $ \eta $. Thus, there are two types of variations with respect to each variable. The variations of the mean variables are:
		\begin{equation*}
		\delta \rho = - \div(\rho \delta \textbf{x}), \qquad \delta \textbf{u} = \dot{(\delta \textbf{x})} - \pd{\textbf{u}}{\textbf{x}} \delta \textbf{x} = \pd{\delta \textbf{x}}{t} + \pd{\delta \textbf{x}}{\textbf{x}} \textbf{u} - \pd{\textbf{u}}{\textbf{x}} \delta \textbf{x}.
		\end{equation*}
		The variation of $ a $ with respect to $ \eta $ is as follows:
		\begin{align*}
		\delta_\eta a &= \int_{t_0}^{t_1} \int_{D(t)} \delta_\eta \left(\rho \frac{\left|\textbf{u}\right|^2}{2} - W\right)dDdt = \int_{t_0}^{t_1} \int_{D(t)} \Bigg( -\pd{W}{\eta} \delta \eta - \pd{W}{\dot{\eta}} \delta \dot{\eta} \Bigg)dDdt\\
		&= \int_{t_0}^{t_1} \int_{D(t)} \Bigg( -\pd{W}{\eta} \delta \eta - \pd{W}{\dot{\eta}} \left( \pd{\delta \eta}{t} + \textbf{u}\cdot \nabla \delta \eta \right) \Bigg)dDdt\\
		&= \int_{t_0}^{t_1} \int_{D(t)} \Bigg( -\pd{W}{\eta} \delta \eta + \pd{}{t}\left(\pd{W}{\dot{\eta}}\right)  \delta \eta + \div\left(\pd{W}{\dot{\eta}} \textbf{u}\right) \delta \eta \Bigg)dDdt\\
		&= \int_{t_0}^{t_1} \int_{D(t)} \Bigg( \left(  -\pd{W}{\eta} + \pd{}{t}\left(\pd{W}{\dot{\eta}}\right)+ \div\left(\pd{W}{\dot{\eta}} \textbf{u}\right) \right) \delta \eta \Bigg)dDdt
		\end{align*}
		Hence, since $ \delta \eta $ vanishes at the boundaries, the Hamilton's principle gives the following equation:
		\begin{equation*}
		\pd{W}{\eta} - \pd{}{t}\left(\pd{W}{\dot{\eta}}\right) - \div\left(\pd{W}{\dot{\eta}} \textbf{u}\right) = 0.
		\end{equation*}
		The variation with respect to $ \dot{\eta} $ is:
		\begin{equation*}
		\delta_\eta \dot{\eta} = \pd{\delta \eta}{t} + \textbf{u}\cdot \nabla \delta \eta.
		\end{equation*}
		The variation of $ \dot{\eta} $ with respect to $ \rho $ is not zero since it is linked to the variation of $ \textbf{u} $:
		\begin{equation*}
		\delta \dot{\eta} = \delta \textbf{u} \cdot \nabla \eta.
		\end{equation*}
		Let us write the Hamilton's principle:
		\begin{align*}
		\delta a &= \int_{t_0}^{t_1} \int_{D(t)} \delta \left(\rho \frac{\left|\textbf{u}\right|^2}{2} - W\right)dDdt = \int_{t_0}^{t_1} \int_{D(t)} \Bigg( \frac{\left|\textbf{u}\right|^2}{2} \delta \rho + \rho \textbf{u} \cdot \delta \textbf{u} - \pd{W}{\rho} \delta \rho - \pd{W}{\dot{\eta}} \delta \dot{\eta} \Bigg)dDdt\\
		&= \int_{t_0}^{t_1} \int_{D(t)} \Bigg( - \frac{\left|\textbf{u}\right|^2}{2} \div (\rho \delta \textbf{x}) + \rho \textbf{u} \cdot \left(\pd{\delta \textbf{x}}{t} + \pd{\delta \textbf{x}}{\textbf{x}} \textbf{u} - \pd{\textbf{u}}{\textbf{x}} \delta \textbf{x}\right) + \pd{W}{\rho} \div (\rho \delta \textbf{x}) - \pd{W}{\dot{\eta}} \delta \dot{\eta} \Bigg)dDdt\\
		& = \int_{t_0}^{t_1} \int_{D(t)} \Bigg( \rho \nabla \frac{\left|\textbf{u}\right|^2}{2} \cdot \delta \textbf{x} + \rho \textbf{u} \cdot \pd{\delta \textbf{x}}{t} - \rho \left(\pd{\textbf{u}}{\textbf{x}}\right)^T \textbf{u} \cdot \delta \textbf{x} + \rho \textbf{u} \cdot \pd{\delta \textbf{x}}{\textbf{x}} \textbf{u} - \rho \nabla \pd{W}{\rho} \cdot \delta \textbf{x} - \pd{W}{\dot{\eta}} \delta \dot{\eta} \Bigg)dDdt\\
		& = \int_{t_0}^{t_1} \int_{D(t)} \Bigg( \rho \textbf{u} \cdot \pd{\delta \textbf{x}}{t} + \rho \textbf{u} \cdot \pd{\delta \textbf{x}}{\textbf{x}} \textbf{u} - \rho \nabla \pd{W}{\rho} \cdot \delta \textbf{x} - \pd{W}{\dot{\eta}} \delta \dot{\eta} \Bigg)dDdt\\
		& = \int_{t_0}^{t_1} \int_{D(t)} \Bigg(  \rho \textbf{u} \cdot \pd{\delta \textbf{x}}{t}  +    \div\big( (\rho \textbf{u} \otimes \delta \textbf{x}) \textbf{u} \big)  -\div (\rho \textbf{u} \otimes \textbf{u}) \cdot \delta \textbf{x} - \rho \nabla \pd{W}{\rho} \cdot \delta \textbf{x} - \pd{W}{\dot{\eta}} \delta \dot{\eta} \Bigg)dDdt\\
		&= \int_{t_0}^{t_1} \int_{D(t)} \Bigg( - \pd{\rho \textbf{u}}{t} \cdot \delta \textbf{x} -\div (\rho \textbf{u} \otimes \textbf{u}) \cdot \delta \textbf{x} - \rho \nabla \pd{W}{\rho} \cdot \delta \textbf{x} - \pd{W}{\dot{\eta}} \delta \dot{\eta} \Bigg)dDdt=\\
		&= \int_{t_0}^{t_1} \int_{D(t)}  - \left( \pd{\rho \textbf{u}}{t} +\div (\rho \textbf{u} \otimes \textbf{u})  + \rho \nabla \pd{W}{\rho} \right) \cdot \delta \textbf{x}~dDdt + \int_{t_0}^{t_1} \int_{D(t)} \Bigg( - \pd{W}{\dot{\eta}} \delta \dot{\eta}  \Bigg)dDdt.
		\end{align*}
		Let us expand the last integral:
		\begin{align*}
		&\int_{t_0}^{t_1} \int_{D(t)} \Bigg( - \pd{W}{\dot{\eta}} \delta \dot{\eta}  \Bigg)dDdt = \int_{t_0}^{t_1} \int_{D(t)} \Bigg( - \pd{W}{\dot{\eta}} \delta \textbf{u} \cdot \nabla \eta \Bigg)dDdt \\
		& =\int_{t_0}^{t_1} \int_{D(t)} \Bigg( - \pd{W}{\dot{\eta}} \left(\pd{\delta \textbf{x}}{t} + \pd{\delta \textbf{x}}{\textbf{x}} \textbf{u} - \pd{\textbf{u}}{\textbf{x}} \delta \textbf{x}\right) \cdot \nabla \eta  \Bigg) dDdt\\
		&= \int_{t_0}^{t_1} \int_{D(t)} \Bigg( - \pd{W}{\dot{\eta}} \left( \pd{\delta \textbf{x}}{t} \cdot \nabla \eta + \div\big( \left( \textbf{u} \otimes \delta \textbf{x} \right) \nabla \eta \big) - \Big( \nabla \eta \div \textbf{u} + \pd{\nabla \eta}{\textbf{x}} \textbf{u} \Big) \cdot \delta \textbf{x} - \left( \pd{\textbf{u}}{\textbf{x}} \right)^T \nabla \eta \cdot \delta \textbf{x} \right) \Bigg)dDdt\\
		&= \int_{t_0}^{t_1} \int_{D(t)} \Bigg(  - \pd{W}{\dot{\eta}} \pd{\delta \textbf{x}}{t} \cdot \nabla \eta - \pd{W}{\dot{\eta}} \div\big( \left( \nabla \eta \cdot \delta\textbf{x}  \right)  \textbf{u} \big) +\pd{W}{\dot{\eta}} \div \textbf{u} \left(\nabla \eta \cdot \delta \textbf{x}\right) \Bigg)dDdt\\
		&\qquad + \int_{t_0}^{t_1} \int_{D(t)} \Bigg( \pd{W}{\dot{\eta}}\pd{\nabla \eta}{\textbf{x}} \textbf{u}  \cdot \delta \textbf{x} + \pd{W}{\dot{\eta}} \left( \pd{\textbf{u}}{\textbf{x}} \right)^T \nabla \eta \cdot \delta \textbf{x} \Bigg)dDdt\\
		&= \int_{t_0}^{t_1} \int_{D(t)} \Bigg(   \pd{}{t} \left(\pd{W}{\dot{\eta}} \nabla \eta\right) \cdot \delta \textbf{x}  + \left( \nabla \eta \cdot \delta\textbf{x}  \right)  \textbf{u} \cdot \nabla \left(\pd{W}{\dot{\eta}}\right)  +\pd{W}{\dot{\eta}} \div \textbf{u} \left(\nabla \eta \cdot \delta \textbf{x}\right)  \Bigg)dDdt\\
		&\qquad + \int_{t_0}^{t_1} \int_{D(t)} \Bigg( \pd{W}{\dot{\eta}}\pd{\nabla \eta}{\textbf{x}} \textbf{u}  \cdot \delta \textbf{x} + \pd{W}{\dot{\eta}} \left( \pd{\textbf{u}}{\textbf{x}} \right)^T \nabla \eta \cdot \delta \textbf{x} \Bigg)dDdt\\
		&= \int_{t_0}^{t_1} \int_{D(t)} \Bigg(   \pd{}{t} \left(\pd{W}{\dot{\eta}} \nabla \eta\right) \cdot \delta \textbf{x}  + \div \left( \pd{W}{\dot{\eta}} \textbf{u} \right) \nabla \eta \cdot \delta \textbf{x} +  \pd{W}{\dot{\eta}}\pd{\nabla \eta}{\textbf{x}} \textbf{u}  \cdot \delta \textbf{x} + \pd{W}{\dot{\eta}} \left( \pd{\textbf{u}}{\textbf{x}} \right)^T \nabla \eta \cdot \delta \textbf{x} \Bigg)dDdt
		\end{align*}
		Thus,
		\begin{align*}
		&\int_{t_0}^{t_1} \int_{D(t)} \Bigg( - \pd{W}{\dot{\eta}} \delta \dot{\eta}  \Bigg)dDdt = \\
		&= \int_{t_0}^{t_1} \int_{D(t)} \Bigg(   \bigg(\pd{}{t} \left(\pd{W}{\dot{\eta}}\right)  + \div \left( \pd{W}{\dot{\eta}} \textbf{u} \right)\bigg) \nabla \eta +  \pd{W}{\dot{\eta}} \bigg( \pd{\nabla \eta}{t}+ \pd{\nabla \eta}{\textbf{x}} \textbf{u}   +  \left( \pd{\textbf{u}}{\textbf{x}} \right)^T \nabla \eta \bigg) \Bigg) \cdot \delta \textbf{x}~dDdt\\
		&= \int_{t_0}^{t_1} \int_{D(t)} \Bigg( \pd{W}{\eta} \nabla \eta + \pd{W}{\dot{\eta}}\dot{(\nabla \eta)}  \Bigg) \cdot \delta \textbf{x}~dDdt = \int_{t_0}^{t_1} \int_{D(t)} \Bigg(\nabla W - \pd{W}{\rho} \nabla \rho  \Bigg)\cdot \delta \textbf{x}~dDdt.
		\end{align*}
		Hence,
		\begin{align*}
		\delta a &= \int_{t_0}^{t_1} \int_{D(t)}  - \left( \pd{\rho \textbf{u}}{t} +\div (\rho \textbf{u} \otimes \textbf{u})  + \rho \nabla \pd{W}{\rho} \right) \cdot \delta \textbf{x}~dDdt + \int_{t_0}^{t_1} \int_{D(t)} \Bigg(\nabla W - \pd{W}{\rho} \nabla \rho  \Bigg)\cdot \delta \textbf{x}~dDdt\\
		&= \int_{t_0}^{t_1} \int_{D(t)}  - \Bigg( \pd{\rho \textbf{u}}{t} +\div (\rho \textbf{u} \otimes \textbf{u})  + \rho \nabla \pd{W}{\rho} + \pd{W}{\rho} \nabla \rho - \nabla W  \Bigg)\cdot \delta \textbf{x}~dDdt\\
		&= \int_{t_0}^{t_1} \int_{D(t)}  - \Bigg( \pd{\rho \textbf{u}}{t} +\div (\rho \textbf{u} \otimes \textbf{u})  + \nabla \left(\rho \pd{W}{\rho} - W \right)  \Bigg)\cdot \delta \textbf{x}~dDdt\\
		&= \int_{t_0}^{t_1} \int_{D(t)}  - \Bigg( \pd{\rho \textbf{u}}{t} +\div \bigg(\rho \textbf{u} \otimes \textbf{u}  + \Big(\rho \pd{W}{\rho} - W \Big) I\bigg)   \Bigg)\cdot \delta \textbf{x}~dDdt.
		\end{align*}
		Eventually, the governing equations read as follows:
		\begin{align*}
		&\pd{\rho}{t} + \div (\rho \textbf{u}) = 0,\\
		&\pd{\rho \textbf{u}}{t} +\div \big(\rho \textbf{u} \otimes \textbf{u}  + p I\big) = 0,\\
		&\pd{}{t}\left(\pd{W}{\dot{\eta}}\right) + \div\left(\pd{W}{\dot{\eta}} \textbf{u}\right) = \pd{W}{\eta},
		\end{align*}
		where the ``pressure'' $ p $ is given by
		\begin{equation*}
		p = \rho \pd{W}{\rho} - W.
		\end{equation*}
		
	\subsection{Hyperbolicity}\label{appendix:extended-model_hyperbolicity}
		In this chapter we suppose that the model is 3-D for the sake of generality, but these results are also valid for the 2-D case. Let us rewrite system \eqref{2-extl_extended-system} in the following form:
		\begin{equation}\label{2-extl_extended-system-primitive-vector}
		\pd{\textbf{U}}{t} + A_x \pd{\textbf{U}}{x} + A_y \pd{\textbf{U}}{y} + A_z \pd{\textbf{U}}{z} = \textbf{S}(\textbf{U}).
		\end{equation}
		where $ \textbf{U} = \left(\rho, u_1, u_2, u_3, \eta, w\right)^T $, and matrices $ A_x $, $ A_y $ and $ A_z $ are:
		\begin{equation*}
		A_x = \left(\begin{matrix}
		u_1                    & \rho & 0 & 0 & 0                       & 0\\
		\frac{p_\rho}{\rho}& u_1 & 0 & 0 & \frac{p_\eta}{\rho} & 0\\
		0                    & 0 & u_1 & 0 & 0                       & 0\\
		0                    & 0 & 0 & u_1 & 0                       & 0\\
		0                    & 0 & 0 & 0 & u_1                       & 0\\
		0                    & 0 & 0 & 0 & 0                       & u_1\\
		\end{matrix}\right), \qquad
		A_y = \left(\begin{matrix}
		u_2                    & \rho & 0 & 0 & 0                       & 0\\
		0                    & u_2 & 0 & 0 & 0                       & 0\\
		\frac{p_\rho}{\rho} & 0 & u_2 & 0 & \frac{p_\eta}{\rho} & 0\\
		0                    & 0 & 0 & u_2 & 0                       & 0\\
		0                    & 0 & 0 & 0 & u_2                       & 0\\
		0                    & 0 & 0 & 0 & 0                       & u_2\\
		\end{matrix}\right),
		\end{equation*}
		\begin{equation*}
		A_z = \left(\begin{matrix}
		u_3                    & \rho & 0 & 0 & 0                       & 0\\
		0                        & u_3 & 0 & 0 & 0                       & 0\\
		0                        & 0 & u_3 & 0 & 0                       & 0\\
		\frac{p_\rho}{\rho} & 0 & 0 & u_3 & \frac{p_\eta}{\rho} & 0\\
		0                        & 0 & 0 & 0 & u_3                       & 0\\
		0                        & 0 & 0 & 0 & 0                       & u_3\\
		\end{matrix}\right).
		\end{equation*}
		Consider a smooth hypersurface $ H(t,x,y,z) = 0 $ and its characteristic vector $ (\xi, \kappa, \chi)^T $ is defined by:
		\begin{equation*}
		\tau = \pd{\rho}{t}, \qquad \xi = \pd{\rho}{x}, \qquad \kappa = \pd{\rho}{y}, \qquad \chi = \pd{\rho}{z}.
		\end{equation*}
		The surface $ H(t,x,y,z) = 0 $ is called \textit{characteristic} if
		\begin{equation*}
		\det \left(\tau I + \xi A_x + \kappa A_y + \chi A_z\right) = 0.
		\end{equation*}
		The system \eqref{2-extl_extended-system-primitive-vector} is $ t $\textit{-hyperbolic} if eigenvalues $ \tau $ of matrix $ \xi A_x + \kappa A_y + \chi A_z $ are real and the corresponding eigenvectors form a basis in $ \mathbb{R}^5 $ \cite{Dafermos2000}. Since the system \eqref{2-extl_extended-system-primitive-vector} is rotationally invariant, one can always transform the unit characteristic vector~$ \frac{\left(\xi, \kappa,\chi\right)^T}{\sqrt{\xi^2 + \kappa^2 + \chi^2}} $ to $ (1,0,0) $. Thus, in order to study the hyperbolicity of a 3-D system it is sufficient to study only the 1-D case, i.e. suppose that $ \textbf{U} = \textbf{U}(x,t) $:
		\begin{equation*}
		\pd{\textbf{U}}{t} + A_x \pd{\textbf{U}}{x} = 0.
		\end{equation*}
		The eigenvalues of $ A_x $ are:
		\begin{equation*}
		\mu_{1,2,3,4} = u, \qquad \mu_{5, 6} = u \pm \sqrt{p_\rho},
		\end{equation*}
		index of $ u_1 $ is omitted here for ease of readability. The corresponding left eigenvectors of are:
		\begin{equation*}
		\begin{aligned}
		\mu = u, \qquad & \textbf{l}_1 = \left( 0, 0, 1, 0, 0, 0 \right)^T,\\
		\mu_2 = u, \qquad & \textbf{l}_2 = \left( 0, 0, 0, 1, 0, 0 \right)^T,\\
		\mu_3 = u, \qquad & \textbf{l}_3 = \left( 0, 0, 0, 0, 1, 0 \right)^T,\\
		\mu_4 = u, \qquad & \textbf{l}_4 = \left( 0, 0, 0, 0, 0, 1 \right)^T,\\
		\mu_5 = u + \sqrt{p_\rho}, \qquad & \textbf{l}_5 = \left( p_\rho, \rho \sqrt{p_\rho} , 0, 0, p_\eta, 0 \right)^T,\\
		\mu_6 = u - \sqrt{p_\rho}, \qquad & \textbf{l}_6 = \left( p_\rho, -\rho \sqrt{p_\rho} , 0, 0, p_\eta, 0 \right)^T.
		\end{aligned}
		\end{equation*}
		The system is hyperbolic, i.e. the eigenvalues are real and the set of $ \textbf{l}_i~(i = 1\dots6) $ is linearly independent if:
		\begin{equation*}
		\pd{p}{\rho} = \rho\big(\rho\varepsilon(\rho)\big)'' + a\lambda \frac{f^2(\eta)}{\rho^2} > 0.
		\end{equation*}
		Let us study the eigenfields of the system. The	right eigenvectors $ \textbf{r}_i~(i = 1\dots6)$ of matrix $ A_x $ are:
		\begin{equation*}
		\begin{aligned}
		\mu = u, \qquad & \textbf{r}_1 = \left( 0, 0, 1, 0, 0, 0 \right)^T,\\
		\mu_2 = u, \qquad & \textbf{r}_2 = \left( 0, 0, 0, 1, 0, 0 \right)^T,\\
		\mu_3 = u, \qquad & \textbf{r}_3 = \left( -p_\eta, 0, 0, 0, p_\rho, 0 \right)^T,\\
		\mu_4 = u, \qquad & \textbf{r}_4 = \left( 0, 0, 0, 0, 0, 1 \right)^T,\\
		\mu_5 = u + \sqrt{p_\rho}, \qquad & \textbf{r}_5 = \left(\rho \sqrt{p_\rho}, p_\rho, 0, 0, 0, 0 \right)^T.\\
		\mu_6 = u - \sqrt{p_\rho}, \qquad & \textbf{r}_6 = \left(-\rho \sqrt{p_\rho}, p_\rho, 0, 0, 0, 0 \right)^T.
		\end{aligned}
		\end{equation*}
		Contact characteristics $ \mu_{1, 2, 3, 4} = u $ are \textit{linearly degenerate}:
		\begin{equation*}
		\mathbf{r}_k\cdot \nabla_{\mathbf{U}} (\mu_k)  \equiv 0, \qquad k = 1,2,3,4, \qquad \nabla_{\textbf{U}} = \left(\partial_\rho, \partial_{u}, \partial_{u_2}, \partial_{u_3}, \partial_\eta, \partial_w\right)^T.
		\end{equation*}
		``Sound'' characteristics $ \mu_{5, 6} = u \pm \sqrt{p_\rho} $ are \textit{genuinely non-linear in the sense of Lax}\cite{Lax2006}:
		\begin{equation*}
		\mathbf{r}_{5, 6}\cdot \nabla_{\mathbf{U}} \mu_{5, 6} = \pm \frac{ 1 }{ 2 } \frac{ \sqrt{p_{\rho}} }{ \rho } \left( 2 + \frac{ \rho p_{\rho\rho} }{ p_{\rho} } \right) \ne 0.
		\end{equation*}

	\subsection{Energy conservation}\label{appendix:extended-model_energy-conservation}
	The equations derived in part \ref{appendix:extended-model-hamilton} of the appendix
	\begin{equation*}
	\begin{aligned}
	&\rho_t + \div(\rho \textbf{u}) = 0,\\
	&(\rho \textbf{u})_t + \div(\rho \textbf{u}\otimes\textbf{u}+pI) = 0,\\
	& \pd{}{t}\left(\pd{W}{\dot{\eta}}\right) + \div \left(\pd{W}{\dot{\eta}}\textbf{u}\right) = \pd{W}{\eta},
	\end{aligned}
	\end{equation*}
	admit the energy conservation law:
	\begin{equation*}
	\left(\rho\frac{|\textbf{u}|^2}{2}+W-\dot{\eta}\pd{W}{\dot{\eta}}\right)_t + \div\left(\rho \frac{|\textbf{u}|^2}{2} \textbf{u} + p\textbf{u} + \Big(W-\dot{\eta}\pd{W}{\dot{\eta}}\Big)\textbf{u}\right) = 0.
	\end{equation*}
	To show this we firstly multiply the second equation by $ \textbf{u} $:
	\begin{equation*}
	\begin{aligned}
	&\textbf{u}\cdot(\rho \textbf{u})_t + \textbf{u} \cdot \div(\rho \textbf{u}\otimes\textbf{u}+pI) = \\
	&\rho\textbf{u}\cdot\textbf{u}_t + \left|\textbf{u}\right|^2\rho_t + \textbf{u}\cdot \left(\textbf{u}\div(\rho \textbf{u})+\pd{\textbf{u}}{\textbf{x}}\rho \textbf{u} \right) + \textbf{u}\cdot\nabla p = 0\\
	& \left(\rho\frac{|\textbf{u}|^2}{2}\right)_t + \underbrace{\frac{|\textbf{u}|^2}{2}\rho_t + \frac{|\textbf{u}|^2}{2}\div(\rho \textbf{u})}_{=0} + \frac{|\textbf{u}|^2}{2}\div(\rho \textbf{u}) + \textbf{u}\cdot\pd{\textbf{u}}{\textbf{x}}\rho \textbf{u} + \textbf{u}\cdot\nabla p = 0\\
	& \left(\rho\frac{|\textbf{u}|^2}{2}\right)_t + \frac{|\textbf{u}|^2}{2}\div(\rho \textbf{u}) +\rho\textbf{u}\cdot \left(\pd{\textbf{u}}{\textbf{x}}\right)^T\textbf{u}+ \ \textbf{u}\cdot\nabla p = 0 \\
	& \left(\rho\frac{|\textbf{u}|^2}{2}\right)_t + \frac{|\textbf{u}|^2}{2}\div(\rho \textbf{u}) +\rho \textbf{u}\cdot\nabla\frac{|\textbf{u}|^2}{2}+ \div(p\textbf{u}) - p\div(\textbf{u}) = 0\\
	& \left(\rho\frac{|\textbf{u}|^2}{2}\right)_t + \div\left(\rho \textbf{u} \frac{|\textbf{u}|^2}{2} + p\textbf{u}\right) - p\div\textbf{u} = 0
	\end{aligned}
	\end{equation*}
	We add the third equation multiplied by $\dot{\eta}$:
	\begin{equation*}
	\begin{aligned}
	&\left(\rho\frac{|\textbf{u}|^2}{2}\right)_t + \div\left(\rho \textbf{u} \frac{|\textbf{u}|^2}{2} + p\textbf{u}\right) +\dot{\eta}\left(\pd{W}{\eta} - \pd{}{t}\left(\pd{W}{\dot{\eta}}\right) - \div \left(\pd{W}{\dot{\eta}}\textbf{u}\right)\right) - p\div\textbf{u} = 0\\
	&\left(\rho\frac{|\textbf{u}|^2}{2}\right)_t + \div\left(\rho \textbf{u} \frac{|\textbf{u}|^2}{2} + p\textbf{u}\right)+\dot{\eta}\pd{W}{\eta} - \dot{\eta}\pd{}{t}\left(\pd{W}{\dot{\eta}}\right) - \div\left(\pd{W}{\dot{\eta}}\textbf{u}\right) - p\div\textbf{u} = 0\\
	&\left(\rho\frac{|\textbf{u}|^2}{2}\right)_t + \div\left(\rho \textbf{u} \frac{|\textbf{u}|^2}{2} + p\textbf{u}\right)+\\
	&\qquad +\dot{\eta}\pd{W}{\dot{\eta}} - \pd{}{t}\left(\dot{\eta}\pd{W}{\dot{\eta}}\right) + \pd{W}{\dot{\eta}}\dot{\eta}_t - \div\left(\dot{\eta}\pd{W}{\dot{\eta}}\textbf{u}\right) + \pd{W}{\dot{\eta}}\textbf{u}\cdot\nabla\dot{\eta} - \rho\pd{W}{\rho}\div\textbf{u} + W\div\textbf{u} = 0\\
	&\left(\rho\frac{|\textbf{u}|^2}{2}+W-\dot{\eta}\pd{W}{\dot{\eta}}\right)_t + \div\left(\rho \textbf{u} \frac{|\textbf{u}|^2}{2} + p\textbf{u} + \left(W-\dot{\eta}\pd{W}{\dot{\eta}}\right)\textbf{u}\right)+\\
	&\qquad+\dot{\eta}\pd{W}{\eta} - \pd{W}{t} + \pd{W}{\dot{\eta}}\dot{\eta}_t - \textbf{u} \cdot \nabla W + \pd{W}{\dot{\eta}}\textbf{u}\cdot \nabla \dot{\eta} - \rho\pd{W}{\rho}\div\textbf{u} = 0
	\end{aligned}
	\end{equation*}
	\begin{equation*}
	\begin{aligned}
	&\left(\rho\frac{|\textbf{u}|^2}{2}+W-\dot{\eta}\pd{W}{\dot{\eta}}\right)_t + \div\left(\rho \textbf{u} \frac{|\textbf{u}|^2}{2} + p\textbf{u} + \left(W-\dot{\eta}\pd{W}{\dot{\eta}}\right)\textbf{u}\right)+\dot{\eta}\pd{W}{\eta}  - \pd{W}{\rho}\rho_t - \pd{W}{\eta}\eta_t\\
	&\qquad - \pd{W}{\dot{\eta}}\dot{\eta}_t + \pd{W}{\dot{\eta}}\dot{\eta}_t  -\pd{W}{\rho} \textbf{u} \cdot \nabla \rho -\pd{W}{\eta} \textbf{u} \cdot \nabla \eta -\pd{W}{\dot{\eta}} \textbf{u} \cdot \nabla \dot{\eta} + \pd{W}{\dot{\eta}}\textbf{u}\cdot \nabla \dot{\eta} - \rho\pd{W}{\rho}\div\textbf{u} = 0\\
	&\left(\rho\frac{|\textbf{u}|^2}{2}+W-\dot{\eta}\pd{W}{\dot{\eta}}\right)_t + \div\left(\rho \textbf{u} \frac{|\textbf{u}|^2}{2} + p\textbf{u} + \left(W-\dot{\eta}\pd{W}{\dot{\eta}}\right)\textbf{u}\right)+\pd{W}{\eta}\dot{\eta} - \pd{W}{\eta}\eta_t\\
	&\qquad -\pd{W}{\eta} \textbf{u} \cdot \nabla \eta  - \pd{W}{\rho}\rho_t -\pd{W}{\rho} \textbf{u} \cdot \nabla \rho - \rho\pd{W}{\rho}\div\textbf{u}  - \pd{W}{\dot{\eta}}\dot{\eta}_t + \pd{W}{\dot{\eta}}\dot{\eta}_t  -\pd{W}{\dot{\eta}} \textbf{u} \cdot \nabla \dot{\eta} + \pd{W}{\dot{\eta}}\textbf{u}\cdot \nabla \dot{\eta}= 0.
	\end{aligned}
	\end{equation*}
	Notice that all the terms except the divergent terms vanish, and only the followig equation is left:
	\begin{equation*}
	\left(\rho\frac{|\textbf{u}|^2}{2}+W-\dot{\eta}\pd{W}{\dot{\eta}}\right)_t + \div\left(\rho \frac{|\textbf{u}|^2}{2} \textbf{u} + p\textbf{u} + \Big(W-\dot{\eta}\pd{W}{\dot{\eta}}\Big)\textbf{u}\right) = 0.
	\end{equation*}